\documentclass[11pt,a4paper]{amsart}
\usepackage{mathpazo}  %Palatino fonts 
\usepackage{amsmath}
\usepackage{amssymb}
\usepackage{graphicx}
\usepackage[all]{xy}
\usepackage{latexsym}
\usepackage{mathrsfs}
\usepackage{fancyhdr}
\usepackage{enumerate}
\usepackage[english]{babel}
\usepackage{lipsum}
\usepackage{etoolbox}
\usepackage[usenames,dvipsnames]{xcolor}

\usepackage{tikz, tkz-euclide, float, pgfplots, overpic}

\usetikzlibrary{arrows,calc, decorations.markings, positioning, fadings}
\usetkzobj{all}
\tikzset{snake it/.style={decorate, decoration=snake}}

\newdimen\tlx
\newdimen\tlx
\newdimen\brx
\newdimen\bry
\usepackage[left=2.5cm, top=3cm, right=2.5cm, bottom=3cm, bindingoffset=0.5cm]{geometry}

\newcommand{\Z}{\ensuremath{\mathbb{Z}}}
\newcommand{\CC}{\ensuremath{\mathcal{C}}} %mathcal C
\DeclareMathOperator{\Hom}{Hom}

\DeclareMathOperator{\End}{End}
\DeclareMathOperator{\add}{add}

\DeclareMathOperator{\ind}{ind} %indecomposables

\newcommand{\modcat}{\operatorname{mod}\nolimits}

\usepackage{collectbox}    % to produce a square box around a number

\makeatletter
\newcommand{\sqbox}{%
    \collectbox{%
        \@tempdima=\dimexpr\width-\totalheight\relax
        \ifdim\@tempdima<\z@
            \fbox{\hbox{\hspace{-.5\@tempdima}\BOXCONTENT\hspace{-.5\@tempdima}}}%
        \else
            \ht\collectedbox=\dimexpr\ht\collectedbox+.5\@tempdima\relax
            \dp\collectedbox=\dimexpr\dp\collectedbox+.5\@tempdima\relax
            \fbox{\BOXCONTENT}%
        \fi
    }%
}
\makeatother

\theoremstyle{theorem}
\newtheorem{theorem}{Theorem}[section]

\theoremstyle{remark}
\newtheorem{remark}[theorem]{Remark}
\newtheorem*{rem}{Remark}
\newtheorem{ex}[theorem]{Example}%[chapter] 
\theoremstyle{definition}
\newtheorem{defi}[theorem]{Definition} %neu

\title{Conway--Coxeter friezes and mutation: a survey}

\author[K.~Baur, E.~Faber, S.~Gratz, K.~Serhiyenko, G.~Todorov]{Karin Baur, Eleonore Faber, Sira Gratz, Khrystyna Serhiyenko, Gordana Todorov}

\address{Institut f\"{u}r Mathematik und Wissenschaftliches Rechnen, 
Universit\"{a}t Graz, NAWI Graz, Heinrichstrasse 36, 
A-8010 Graz, Austria}
\email{baurk@uni-graz.at}

\address{
School of Mathematics, University of Leeds, Leeds, LS2 9EJ, UK
}
\email{e.m.faber@leeds.ac.uk}

\address{
Mathematical Institute, University of Oxford, Oxford, OX2 6GG, UK
}
\email{gratz@maths.ox.ac.uk}

\address{Department of Mathematics, University of California at Berkeley, Berkeley, CA 94720, USA}
\email{khrystyna.serhiyenko@berkeley.edu}

\address{Department of Mathematics, Northeastern University, Boston, MA 02115, USA}
\email{g.todorov@neu.edu}
\date{\today}
\begin{document}

\begin{abstract}
In this survey article we explain the intricate links between Conway-Coxeter friezes and cluster 
combinatorics. 
More precisely, we provide a formula, relying solely on the shape of the frieze, 
describing how each individual entry in the frieze changes under cluster mutation. 
Moreover, we provide a combinatorial formula for the number of submodules of a string module, 
and with that a simple way to compute the frieze associated to a fixed cluster tilting object in a 
cluster category of Dynkin type $A$ in the sense of Caldero and Chapoton. \\

\end{abstract}

\maketitle

{\footnotesize
\tableofcontents 
}

%%%%%%%%%%%%
\section{Introduction}\label{sec:introduction}

Cluster algebras were introduced by Fomin and Zelevinsky in \cite{FZ1}. A key motivation was to provide 
an algebraic framework for phenomena observed in the study of dual canonical bases for quantised 
enveloping algebras and in total positivity for reductive groups. 

Cluster categories were introduced in 2005, \cite{BMRRT}, \cite{CalderoChapotonSchiffler} 
to give a categorical 
interpretation of cluster algebras. The following table shows the beautiful interplay and 
correspondences between 
cluster algebras and cluster categories in type $A$. Note that the correspondences between the 
first and second column hold more generally, not only in type $A$:  
Caldero and Chapoton \cite{CalderoChapoton} have provided a formal link between cluster 
categories and cluster algebras by introducing what is now most commonly known as the {\em Caldero Chapoton map} (short: CC-map) or {\em cluster character}. 
%, denoted by $\rho$ in the table. 
Fixing a cluster tilting object (which takes on the role of the initial cluster), it associates to each indecomposable in the cluster category a unique cluster variable in the associated cluster algebra, sending the indecomposable summands of the cluster tilting object to the initial cluster. 

\begin{center}
\begin{tabular}{|c|c|c|c|}
\hline 
Cluster algebra & $\leftarrow$ & Cluster category &  Polygon \\
\hline 
cluster variables  & CC-map & indecomposable objects  & diagonals \\
clusters && cluster tilting objects & triangulations \\
mutations && mutations & flip \\
\hline
\end{tabular}
\end{center}

In the 70s, Coxeter and Conway first studied frieze patterns of numbers (\cite{CoCo1} and \cite{CoCo2}). When these numbers 
are positive integers, they showed that the frieze patterns 
arise from triangulations of polygons. Thus we can extend this table by a further column: 
\begin{center}
\begin{tabular}{|c|c|c|}
\hline
$\dots$ &  Polygon & Frieze\\
\hline 
  & diagonals & integers \\
 & triangulations & sequences of 1's \\
 & flip & ?? \\
\hline
\end{tabular}
\end{center}

Here the last entry is missing: the meaning of mutation or flip on the level of frieze patterns was 
not known until now. 
The purpose of this survey article is to show how to 
complete the picture of cluster combinatorics in the context of friezes. It is based on the 
paper~\cite{BFGST17} where more background on cluster categories can be found and where 
the proofs are included. 

More precisely, we determine how mutation of a cluster affects the associated frieze, thus 
effectively introducing the notion of a mutation of friezes that is compatible with mutation in the 
associated cluster algebra. 
This provides a useful new tool to study cluster combinatorics of Dynkin type $A$.

In order to deal with the mutations for friezes we will use cluster categories and generalized cluster categories as introduced by Buan, Marsh, Reineke, Reiten and Todorov \cite{BMRRT} for hereditary algebras and by Amiot \cite{[Am]}  more generally. In both cases, cluster categories are triangulated categories in which the combinatorics of cluster algebras receives a categorical interpretation: cluster variables correspond to rigid 
indecomposable objects  
and clusters correspond to cluster tilting objects. One of the essential features in the definition of cluster algebras is the process of mutation, which replaces one element of the cluster by a another unique  element  such that a new cluster is created. The corresponding categorical mutation replaces an indecomposable summand of a cluster tilting object by another unique indecomposable object using approximations in the triangulated categories; this process creates another cluster tilting object which corresponds to the mutated cluster.

We now explain the different players appearing in the table above. 

%%%%%%%%%%%%
\subsection{Frieze patterns}\label{ssec:fr-pattern}

The notion of friezes was introduced by Coxeter \cite{Coxeter}; it was Gauss's \emph{pentagramma mirificum} which was the original inspiration. We recall that a \emph{frieze} is a grid of positive integers consisting of a finite number of infinite rows: the top and bottom rows are infinite rows of $0$s and the second to top and second to bottom are infinite rows of $1$s as one can see on the following diagram

\[
 \xymatrix@=0.5em{ 
 &\ldots && 0 && 0 && 0 && 0 && \ldots &\\
   && 1 && 1 && 1 && 1 && 1 && \\
  &\ldots && m_{-1,-1} && m_{00} && m_{11} && m_{22} && \ldots& \\
   && m_{-2,-1} && m_{-1,0} && m_{01} && m_{12} && m_{23} && \\
   & \ldots && \ldots && \ldots && \ldots && \ldots && \ldots \\
    && 1 && 1 && 1 && 1 && 1 && \\
   &\ldots && 0 && 0 && 0 && 0 && \ldots &
  }
\]
The entries of the frieze satisfy the {\em frieze rule}: for every set of adjacent numbers arranged in a {\em diamond}
\[
 \xymatrix@=0.5em{
 & b & \\
 a && d \\
 & c &
 }
\]
we have
\[
 ad - bc = 1.
\]

The sequence of integers  in the first non-trivial row, $(m_{ii})_{i \in \Z}$, is called 
{\em quiddity sequence}. This sequence completely determines the frieze. Each frieze is also periodic, since it is invariant under glide reflection. The order of the frieze is defined to be the number of rows minus one. It follows that each frieze of order $n$ is $n$-periodic.

Among the famous results about friezes is the bijection between the friezes of order $n$ and triangulations of a convex $n$-gon, which was proved by Conway and Coxeter in \cite{CoCo1} and \cite{CoCo2}. This was used to set the first link with cluster combinatorics using \cite{CalderoChapotonSchiffler} and \cite{CalderoChapoton} by Caldero and Chapoton. 
Recently, frieze patterns have been generalized in several directions and found applications in various areas of mathematics, for an overview see \cite{Morier-Genoud}.
%\sg{Sounds good.}

%%%%%%%%%%%%
\subsection{Cluster algebras}

Fomin and Zelevinsky introduced the notion of cluster algebras in \cite{FZ1}. Cluster algebras are commutative algebras generated by \emph{cluster variables}; cluster variables are obtained from an \emph{initial cluster} (of variables) by replacing one element at a time according to a prescribed rule, where the rule is given either by a skew-symmetric (or more generally skew-symmetrizable) matrix or, equivalently by a quiver with no loops nor $2$-cycles. The process of replacing one element of a cluster by another unique element in order to obtain another cluster, together with the prescribed change of the quiver, is called \emph{mutation}. Finite sequences of iterated mutations create new clusters and new cluster variables; all cluster variables are obtained in such a way. 

The process of such mutations may never stop, however if the quiver is of Dynkin type, then by the theorem of Fomin and Zelevinsky, this process stops and one obtains a finite number of cluster variables \cite{FZ2}. Among those cluster algebras, the best behaved and understood are the cluster algebras of type $A$. The clusters of the cluster algebra of type $A_{n-3}$ are in bijection with the triangulations of a convex $n$-gon, for $n\geq 3$. This is exactly what is employed in this work in order to relate and use cluster categories, via triangulations of $n$-gon, so that we can describe the mutations of friezes of order $n$.
Since we will also be dealing  with quivers $Q'$ which are mutation equivalent to the quivers of type 
$A_n$ and may have nontrivial potential, we need to consider 
generalized cluster categories $\mathcal C_{(Q',W)}$, which are shown to be triangle equivalent 
to $\mathcal C_Q$ \cite{[Am]}.

%%%%%%%%%%%%
\subsection{Cluster categories}

Let $Q$ be an acyclic quiver with $n$ vertices, over an algebraically closed field. We consider 
the category $\modcat kQ$ of (finitely generated) modules over $kQ$, or, equivalently, the 
category rep$Q$ of representations of the quiver $Q$. The bounded derived category 
$D^b(kQ)$ can be viewed as $\cup_{i\in\mathbb{Z}} \modcat(kQ)[i]$, with connecting morphisms. 

As an example, consider the quiver 
\[
\xymatrix{Q: &1 & 2 \ar[l]& 3\ar[l]
}
\]
The module category of the path algebra $kQ$ has six indecomposable objects up to isomorphisms, 
with irreducible maps between them as follows: 
\[
\xymatrix@!C=0pt@!R=0pt{
 && P_3\ar[rd] \\
& P_2 \ar[rd]\ar[ru] && I_2\ar[rd] \\
P_1\ar[ru] && S_2\ar[ru] && I_3
}
\]
The modules $P_i$ are indecomposable projective, the $I_i$ are indecomposable injectives 
and the $S_i$ are the simple modules, with $I_1=P_3$, $S_1=P_1$ and $S_3=I_3$. The bounded 
derived category then looks as follows (the arrows indicate the connecting morphisms): 
\[
\includegraphics[width=12cm]{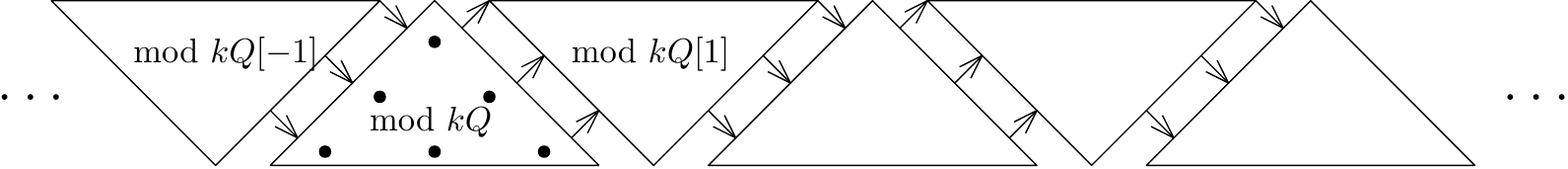}
\]

Let $Q$ be a Dynkin quiver of type $A$. Let $\CC$ be the associate cluster category, which by definition is $\CC = \CC_Q = D^b(kQ)/\tau^{-1}[1]$ where  $D^b(kQ)$ is the bounded derived category of the path algebra $kQ$ with the suspension functor $[1]$ and the Auslander-Reiten functor $\tau$. In this case, the specialized CC map gives a direct connection between the Auslander-Reiten quiver of the cluster category $\CC$ with a fixed cluster tilting object $T$, and the associated frieze $F(T)$ in the following way: recall that each vertex of the Auslander-Reiten quiver corresponds to an isomorphism class of an indecomposable object in the cluster category. When the specialized CC map is applied to a representative of each isomorphism class and the vertex is labeled by that value, one only needs to complete those rows by the rows of $1$s and $0$s at the top and bottom in order to obtain a frieze, cf.\ \cite[Proposition~5.2]{CalderoChapoton}.

%%%%%%%%%%%%

\section{From cluster categories to frieze patterns}

Let $\CC$ be a cluster category, let $T$ be a cluster tilting object and let  $B_T = \End_{\CC}(T)$ be the endomorphism algebra, which is also called a \emph{cluster-tilted algebra}. The module category $\modcat (B_T)$ is shown to be equivalent to the quotient category $\CC \big/ \add(T[1])$ of the cluster category. This result by Buan, Marsh and Reiten is used in a very essential way: each indecomposable object in $\CC$, which is not isomorphic to a summand of $T[1]$ corresponds to an indecomposable $B_T$-module, preserving the structure of the corresponding Auslander-Reiten quivers; at the same time the indecomposable summands of $T[1]$ correspond to the suspensions of the indecomposable projective $B_T$-modules in the generalized cluster category of the algebra $B_T$. 

When $\CC$ is the cluster category associated to the Dynkin quiver of type $A$, for each cluster tilting object $T$, the associated specialized CC-map sends each indecomposable summand of $T[1]$ to $1$ and each indecomposable $B_T$-module $M$ to the number of its submodules, as we explain now.
In the actual Caldero-Chapoton formula for cluster variable $x_M$ in terms of the initial cluster variables, the coefficients are given as the Euler-Poincar{\'e} characteristics of the Grassmannians of submodules of the module $M$. In this expression the sum is being taken over the dimension vectors of the submodules of $M$. However in this set-up, since all indecomposable $B_T$-modules are string modules, all the Grassmannians are just points. 
The {\em specialized Caldero-Chapoton map} is the map we get from postcomposing the 
CC-map associated to $T$ with the specialization of the initial cluster variables to one. 
Hence the sum is equal to the number of submodules and the values of the specialized CC-map are positive integers. 
The values of the specialized CC-map are now entered in the AR-quiver of the cluster category $\CC$ at the places of the corresponding indecomposable objects.
The image of this generalized CC-map only needs to be completed with the rows of $1$s and $0$s above and below in order to obtain the {\em frieze associated to the cluster tilting object} $T$, denoted by $F(T)$.

\begin{figure}
%If this picture is not centered, probably it is because it is too wide. Tune the scale parameter in the next line.
\begin{tikzpicture}[scale=0.85,line join = round, line cap = round,>=stealth]
%Define parameters for scalability
\def\hh{3}
\def\xone{1.7}
\def\nbhd{0.15}
\def\nbhdMore{0.188} %unfortunately needs to be handtuned...
\def\hratio{0.8}
\def\xratio{1}
\def\endgapmultiplier{4.2}
%pgf macros to get scalable coordinates for other points. The final points will be centered.
\pgfmathsetmacro\cxSLpOne{\xratio*\xone};
\pgfmathsetmacro\cxSLmOne{-\xratio*\xone};

\pgfmathsetmacro\cxblankright{2*\xratio*\xone};
\pgfmathsetmacro\cyblankright{\hratio*\hh};
\pgfmathsetmacro\cxblankleft{-2*\xratio*\xone};
\pgfmathsetmacro\cyblankleft{\hratio*\hh};
\pgfmathsetmacro\heightend{\hh*0.5};
\pgfmathsetmacro\rightend{\endgapmultiplier*\xratio*\xone};
\pgfmathsetmacro\leftend{-\endgapmultiplier*\xratio*\xone};
\pgfmathsetmacro\wiggledist{\rightend-\cxblankright};
\pgfmathsetmacro\wigglestep{\wiggledist*0.125};
%Coordinates depend on parameters. The picture is bounded below by y=0.
%The picture is centered at x=0.
\coordinate (SL) at (0, \hh);
\coordinate (SLp1) at (\cxSLpOne, 0);
\coordinate (SLm1) at (\cxSLmOne, 0);
\coordinate (Sblankright) at (\cxblankright, \cyblankright);
\coordinate (Sblankleft) at (\cxblankleft, \cyblankleft);
\coordinate (Srightend) at (\rightend,\heightend);
\coordinate (Sleftend) at (\leftend,\heightend);
%small dots
\coordinate (nbhdSLleft) at ($(SL) !\nbhd! (SLm1)$);
\coordinate (nbhdSLright) at ($(SL) !\nbhd! (SLp1)$);
\coordinate (nbhdSLm1right) at ($(SLm1) !\nbhd! (SL)$);
\coordinate (nbhdSLm1left) at ($(SLm1) !\nbhdMore! (Sblankleft)$);
\coordinate (nbhdSLp1left) at ($(SLp1) !\nbhd! (SL)$);
\coordinate (nbhdSLp1right) at ($(SLp1) !\nbhdMore! (Sblankright)$);

%encircle
%left
\coordinate (leftc1) at ($(nbhdSLleft)!0.4!(SL)$);
\coordinate (leftc2) at ($(nbhdSLleft)!0.3!(nbhdSLright)$);
\coordinate (leftc3) at ($(SLm1)!0.5!(SL) + (0.3,0)$);
\coordinate (leftc5help) at ($(nbhdSLm1right)!0.5!(SLm1)$);
\coordinate (leftc5) at ($(leftc5help)!1.8!(SLm1)$);
\coordinate (leftc4help) at ($(leftc5help)+(0.3,0)$);
\coordinate (leftc4) at ($(leftc3)!1.2!(leftc4help)$);
\coordinate (leftc6help) at ($(nbhdSLm1right)!0.5!(nbhdSLm1left)$);
\coordinate (leftc7) at ($(SLm1)!0.5!(SL) - (0.3,0)$);
\coordinate (leftc6) at ($(leftc7)!1.2!(leftc6help)$);
\coordinate (leftc8) at ($(leftc1)- (0.25,0)$);

%right
\coordinate (rightc1) at ($(nbhdSLright)!0.5!(SL)$);
\coordinate (rightc2) at ($(nbhdSLright)!0.3!(nbhdSLleft)$);
\coordinate (rightc3) at ($(SLp1)!0.5!(SL) - (0.3,0)$);
\coordinate (rightc5help) at ($(nbhdSLp1left)!0.5!(SLp1)$);
\coordinate (rightc5) at ($(rightc5help)!1.8!(SLp1)$);
\coordinate (rightc4help) at ($(rightc5help)-(0.3,0)$);
\coordinate (rightc4) at ($(rightc3)!1.2!(rightc4help)$);
\coordinate (rightc6help) at ($(nbhdSLp1right)!0.5!(nbhdSLp1left)$);
\coordinate (rightc7) at ($(SLp1)!0.5!(SL) + (0.3,0)$);
\coordinate (rightc6) at ($(rightc7)!1.2!(rightc6help)$);
\coordinate (rightc8) at ($(rightc1) + (0.25,0)$);

%Nodes and labels
\node[above] at (SL) {$S_{\ell}$};
\node[below] at (SLp1) {$S_{\ell+1}$};
\node[below] at (SLm1) {$S_{\ell-1}$};
%\node at (Sblankright) {hello};
%\node at (Sblankleft) {what?};
\node[right] at (Srightend) {$S_{m+1}$};
\node[left] at (Sleftend) {$S_{1}$};

%Drawings
\draw[Blue,line width=1.6] (SL)--(SLp1);
\draw[Blue,line width=1.6] (SL)--(SLm1);
\draw[Blue,line width=1.6] (SLm1)--(Sblankleft);
\draw[Blue,line width=1.6] (SLp1)--(Sblankright);

\draw[Blue, opacity=1,line width =1.6,path fading=east] plot[smooth] coordinates {
    (Sblankright) 
    ($(Srightend) + (-\wiggledist +\wigglestep, -0.21*\hh)$)
     ($(Srightend) + (-\wiggledist +2*\wigglestep, +0.18*\hh)$)
     ($(Srightend) + (-\wiggledist +3*\wigglestep, -0.15*\hh)$)
     ($(Srightend) + (-\wiggledist +4*\wigglestep, +0.12*\hh)$)
     ($(Srightend) + (-\wiggledist +5*\wigglestep, -0.09*\hh)$)
     ($(Srightend) + (-\wiggledist +6*\wigglestep, +0.06*\hh)$)
     ($(Srightend) + (-\wiggledist +7*\wigglestep,  -0.03*\hh)$)
     (Srightend)
   }; 
\draw[Blue, opacity=1,line width =1.6,path fading=west] plot[smooth] coordinates {
    (Sblankleft) 
    ($(Sleftend) + (\wiggledist -\wigglestep, -0.21*\hh)$)
     ($(Sleftend) + (\wiggledist -2*\wigglestep, +0.18*\hh)$)
     ($(Sleftend) + (\wiggledist -3*\wigglestep, -0.15*\hh)$)
     ($(Sleftend) + (\wiggledist -4*\wigglestep, +0.12*\hh)$)
     ($(Sleftend) + (\wiggledist -5*\wigglestep, -0.09*\hh)$)
     ($(Sleftend) + (\wiggledist -6*\wigglestep, +0.06*\hh)$)
     ($(Sleftend) + (\wiggledist -7*\wigglestep,  -0.03*\hh)$)
     (Sleftend)
   }; 
   
%Enclosing paths
\draw[Red, opacity=1,line width =1.4, dashed] plot[smooth] coordinates {
    (leftc1) 
    (leftc2)
    (leftc3)
    (leftc4)
    (leftc5)
    (leftc6)
    (leftc7)
    (leftc8)
    (leftc1)
   }; 
\node[right] at ($(leftc2)!0.8!(leftc4)$) {\color{Red} $N_{\ell-1}$};
\draw[Red, opacity=1,line width =1.4, dashed] plot[smooth] coordinates {
    (rightc1) 
    (rightc2)
    (rightc3)
    (rightc4)
    (rightc5)
    (rightc6)
    (rightc7)
    (rightc8)
    (rightc1)
   }; 
\node[left] at ($(rightc2)!0.8!(rightc4)$) {\color{Red} $N_{\ell}$};
%Marks
\foreach \Point in {(SL), (SLp1), (SLm1), (Sblankright), (Sblankleft)}{
\draw[LimeGreen] \Point circle[radius=3pt];
\fill[LimeGreen] \Point circle[radius=3pt];
}
\foreach \Point in {(Srightend),(Sleftend)}{
\draw[opacity=0.2,LimeGreen] \Point circle[radius=3pt];
\fill[opacity=0.2,LimeGreen] \Point circle[radius=3pt];
}
\foreach \Point in {(nbhdSLright),(nbhdSLp1left),(nbhdSLleft),(nbhdSLm1right),(nbhdSLp1right),($(Sblankright) !\nbhd! (SLp1)$),(nbhdSLm1left),($(Sblankleft) !\nbhd! (SLm1)$)}{
\draw[Green] \Point circle[radius=2.4pt];
\fill[Green] \Point circle[radius=2.4pt];
}
%TestMarks
%\foreach \Point in {(leftc1),(leftc2),(leftc3),(leftc4),(leftc5),(leftc6),(leftc7),(leftc8)}{
%\draw \Point circle[radius=1pt];
%}
\end{tikzpicture}
\caption{A string module $M=(k_1, \ldots, k_m)$ with legs $N_l$.}
\label{F:string}
\end{figure}
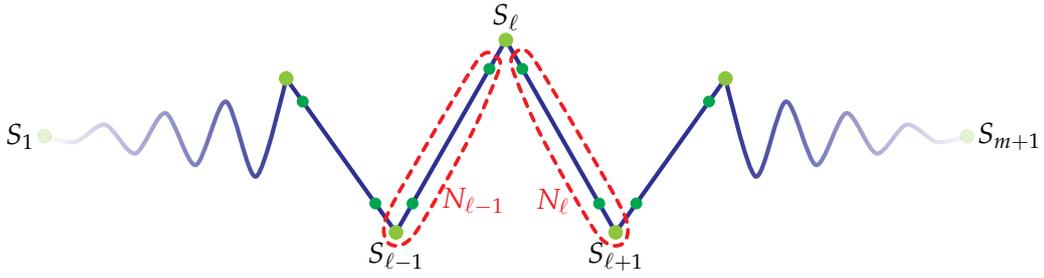

Since the generalized CC-map for cluster categories of Dynkin type $A$ is given in terms of the number of submodules of $B_T$-modules, the first goal of the paper is to give a formula for the number of submodules.
This is determined by the following result, hence providing a combinatorial formula for the number of submodules of any given indecomposable $B_T$-module. Its proof can be 
found in~\cite[Section 4]{BFGST17}. 
We recall that 
each $B_T$-module is a string module and hence has a description in terms of 
the lengths of the individual legs. 
Let $(k_1,\dots, k_m)$ denote these lengths, cf.\ Figure~\ref{F:string}. We further denote by $s(M)$ the number of submodules of a $B_T$-module $M$.

%%%%%%% String module - see figure-sketch.tex
%\begin{figure}[h]
%\includegraphics[width=\textwidth]{string1.png}  %The quality is not very good, but I don't have photoshop, so the pdf does not work - I include the source code blow
%\caption{A string module $M$ of shape $(k_1, \ldots, k_m)$ with legs $N_l$.}
%\label{F:string}
%\end{figure}

%%%%%% String module - see figure-sketch.tex
%%%%%%%%%%%%%%%%%%%%%%%%

\begin{theorem} 
Let $M$ be an indecomposable $B_T$-module, of shape $(k_1,\dots,k_m)$. Then 
the number of submodules of $M$ is given as: 
\[
 s(M) = 1 + \sum_{j=0}^m\sum_{|I|=m-j}\prod_{i\in I}k_i \ ,
\]
where the second sum runs over all admissible subsets $I$ of $\{1,\dots, m\}$. 
\end{theorem}

The above formula for the number of submodules of an indecomposable string module of the shape $(k_1,\ldots,k_m)$ is given in terms of those integers, which makes it clear that 
we need to determine those integers.
And that is exactly what was done: using the position of the module in the AR-quiver and the information about the positions of the indecomposable projective $B_T$-modules, the procedure for finding the numerical invariants $(k_1,\ldots,k_m)$ of the module was given. This purely combinatorial way of computing the numbers of submodules was the basis for computing the associated friezes, and eventually, mutations of friezes.

We end this section by giving an example illustrating the frieze pattern obtained through the 
specialized CC-map. 
\begin{ex}\label{E:category_frieze}

We now illustrate several notions on the example of the cluster category $\CC_{A_{11}}$: 
the Auslander-Reiten quiver of $\CC_{A_{11}}$, a cluster tilting object $T$, 
the cluster-tilted algebra $B_T$ and the Auslander-Reiten quiver of the 
generalized cluster category of $B_T$ where the modules are given by their composition 
factors. 
The Auslander-Reiten quiver of $\CC_{A_{11}}$ is the quotient of 
the Auslander-Reiten quiver of $D^b(kA_{11})$ by the action of $\tau^{-1}[1]$, 
a fundamental domain for which is depicted in black below. 
We pick the cluster tilting object $T = \bigoplus_{i=1}^{11} T_i$ whose indecomposable 
summands are marked with circles:

\vskip-2ex
\begin{center}
\begin{tikzpicture}[scale=0.6, transform shape]

	%first row
	\node[lightgray] (Z1) at (0,10) {$\bullet$};
	\node[lightgray] (Z2) at (2,10) {$\bullet$};
	\node[lightgray, circle, draw] (Z3) at (4,10) {$\bullet$};
	\node[lightgray] (Z4) at (6,10) {$\bullet$};
	\node[lightgray] (Z5) at (8,10) {$\bullet$};
	\node[circle, draw] (Z6) at (10,10) {$\bullet$};
	\node (Z7) at (12,10) {$\bullet$};
	\node[lightgray, circle, draw] (Z8) at (14,10) {$\bullet$};
	\node[lightgray] (Z9) at (16,10) {$\bullet$};
	\node[lightgray, circle, draw] (Z10) at (18,10) {$\bullet$};
	\node[lightgray] (Z11) at (20,10) {$\bullet$};
	\node[lightgray] (Z12) at (22,10) {$\bullet$};
	
	%2nd row
	\node[lightgray] (Y1) at (1,9) {$\bullet$};
	\node[lightgray] (Y2) at (3,9) {$\bullet$};
	\node[lightgray, circle, draw] (Y3) at (5,9) {$\bullet$};
	\node[lightgray] (Y4) at (7,9) {$\bullet$};
	\node (Y5) at (9,9) {$\bullet$};
	\node[] (Y6) at (11,9) {$\bullet$};
	\node (Y7) at (13,9) {$\bullet$};
	\node[lightgray] (Y8) at (15,9) {$\bullet$};
	\node[lightgray] (Y9) at (17,9) {$\bullet$}; 
	\node[lightgray, circle, draw] (Y10) at (19,9) {$\bullet$};
	\node[lightgray] (Y11) at (21,9) {$\bullet$};
	%\node[lightgray] (Y12) at (23,9) {$\bullet$};
	
	%3rd row
	\node[lightgray] (S1) at (0,8) {$\bullet$};
	\node[lightgray] (S2) at (2,8) {$\bullet$};
	\node[lightgray] (S3) at (4,8) {$\bullet$};
	\node[lightgray, ] (S4) at (6,8) {$\bullet$};
	\node[] (S5) at (8,8) {$\bullet$};
	\node (S6) at (10,8) {$\bullet$};
	\node (S7) at (12,8) {$\bullet$};
	\node[] (S8) at (14,8) {$\bullet$};
	\node[lightgray, ] (S9) at (16,8) {$\bullet$};
	\node[lightgray] (S10) at (18,8) {$\bullet$};
	\node[lightgray] (S11) at (20,8) {$\bullet$};
	\node[lightgray] (S12) at (22,8) {$\bullet$};
	
	%4th row
	\node[lightgray] (B1) at (1,7) {$\bullet$};
	\node[lightgray] (B2) at (3,7) {$\bullet$};
	\node[lightgray] (B3) at (5,7) {$\bullet$};
	\node[circle, draw] (B4) at (7,7) {$\bullet$};
	\node (B5) at (9,7) {$\bullet$};
	\node (B6) at (11,7) {$\bullet$};
	\node (B7) at (13,7) {$\bullet$};
	\node (B8) at (15,7) {$\bullet$};
	\node[lightgray, circle, draw] (B9) at (17,7) {$\bullet$}; 
	\node[lightgray] (B10) at (19,7) {$\bullet$};
	\node[lightgray] (B11) at (21,7) {$\bullet$};
	%\node[lightgray] (B12) at (23,7) {$\bullet$};
	
	%5th row
	\node[lightgray] (C1) at (0,6) {$\bullet$};
	\node[lightgray] (C2) at (2,6) {$\bullet$};
	\node[lightgray] (C3) at (4,6) {$\bullet$};
	\node (C4) at (6,6) {$\bullet$};
	\node (C5) at (8,6) {$\bullet$};
	\node (C6) at (10,6) {$\bullet$};
	\node (C7) at (12,6) {$\bullet$};
	\node (C8) at (14,6) {$\bullet$};
	\node[] (C9) at (16,6) {$\bullet$};
	\node[lightgray] (C10) at (18,6) {$\bullet$};
	\node[lightgray] (C11) at (20,6) {$\bullet$};
	\node[lightgray] (C12) at (22,6) {$\bullet$};
	
	%6th row
	\node[lightgray] (D1) at (1,5) {$\bullet$};
	\node[lightgray] (D2) at (3,5) {$\bullet$};
	\node (D3) at (5,5) {$\bullet$};
	\node[] (D4) at (7,5) {$\bullet$};
	\node (D5) at (9,5) {$\bullet$};
	\node (D6) at (11,5) {$\bullet$};
	\node (D7) at (13,5) {$\bullet$};
	\node (D8) at (15,5) {$\bullet$};
	\node (D9) at (17,5) {$\bullet$}; 
	\node[lightgray] (D10) at (19,5) {$\bullet$};
	\node[lightgray] (D11) at (21,5) {$\bullet$};
	%\node[lightgray] (D12) at (23,5) {$\bullet$};
	
	%7th row
	\node[lightgray] (E1) at (0,4) {$\bullet$};
	\node[lightgray] (E2) at (2,4) {$\bullet$};
	\node (E3) at (4,4) {$\bullet$};
	\node (E4) at (6,4) {$\bullet$};
	\node (E5) at (8,4) {$\bullet$};
	\node (E6) at (10,4) {$\bullet$};
	\node (E7) at (12,4) {$\bullet$};
	\node (E8) at (14,4) {$\bullet$};
	\node (E9) at (16,4) {$\bullet$};
	\node (E10) at (18,4) {$\bullet$};
	\node[lightgray] (E11) at (20,4) {$\bullet$};
	\node[lightgray] (E12) at (22,4) {$\bullet$};
	
	%8th row
	\node[lightgray] (F1) at (1,3) {$\bullet$};
	\node[circle, draw] (F2) at (3,3) {$\bullet$};
	\node[] (F3) at (5,3) {$\bullet$};
	\node (F4) at (7,3) {$\bullet$};
	\node (F5) at (9,3) {$\bullet$};
	\node (F6) at (11,3) {$\bullet$};
	\node (F7) at (13,3) {$\bullet$};
	\node (F8) at (15,3) {$\bullet$};
	\node[] (F9) at (17,3) {$\bullet$};
	\node (F10) at (19,3) {$\bullet$};
	\node[lightgray, circle, draw] (F11) at (21,3) {$\bullet$};
	%\node[lightgray, ] (F12) at (23,3) {$\bullet$}; 
	
	%9th row
	\node[lightgray] (G1) at (0,2) {$\bullet$};
	\node (G2) at (2,2) {$\bullet$};
	\node (G3) at (4,2) {$\bullet$};
	\node (G4) at (6,2) {$\bullet$};
	\node (G5) at (8,2) {$\bullet$};
	\node[] (G6) at (10,2) {$\bullet$};
	\node[circle, draw] (G7) at (12,2) {$\bullet$};
	\node (G8) at (14,2) {$\bullet$};
	\node (G9) at (16,2) {$\bullet$};
	\node (G10) at (18,2) {$\bullet$};
	\node (G11) at (20,2) {$\bullet$};
	\node[lightgray] (G12) at (22,2) {$\bullet$};
	
	%10th row
	\node (H1) at (1,1) {$\bullet$};
	\node[] (H2) at (3,1) {$\bullet$};
	\node[circle, draw] (H3) at (5,1) {$\bullet$};
	\node[] (H4) at (7,1) {$\bullet$};
	\node (H5) at (9,1) {$\bullet$};
	\node[] (H6) at (11,1) {$\bullet$};
	\node (H7) at (13,1) {$\bullet$};
	\node[] (H8) at (15,1) {$\bullet$};
	\node (H9) at (17,1) {$\bullet$}; 
	\node[circle, draw] (H10) at (19,1) {$\bullet$};
	\node[] (H11) at (21,1) {$\bullet$};
	%\node[lightgray] (H12) at (23,1) {$\bullet$};
	
	%11th row
	\node[circle, draw] (I1) at (0,0) {$\bullet$};
	\node[] (I2) at (2,0) {$\bullet$};
	\node[circle, draw] (I3) at (4,0) {$\bullet$};
	\node[] (I4) at (6,0) {$\bullet$};
	\node[] (I5) at (8,0) {$\bullet$};
	\node[circle, draw] (I6) at (10,0) {$\bullet$};
	\node[] (I7) at (12,0) {$\bullet$};
	\node[circle, draw] (I8) at (14,0) {$\bullet$};
	\node[]  (I9) at (16,0) {$\bullet$}; 
	\node[circle, draw] (I10) at (18,0) {$\bullet$};
	\node (I11) at (20,0) {$\bullet$};
	\node[] (I12) at (22,0) {$\bullet$};
	
	%arrows 1st to 2nd
	\draw[lightgray, ->] (Z1) -- (Y1);
	\draw[lightgray,->] (Z2) -- (Y2);
	\draw[lightgray,->] (Z3) -- (Y3);
	\draw[lightgray,->] (Z4) -- (Y4);
	\draw[lightgray,->] (Z5) -- (Y5);
	\draw[->] (Z6) -- (Y6);
	\draw[->] (Z7) -- (Y7);
	\draw[lightgray,->] (Z8) -- (Y8);
	\draw[lightgray,->] (Z9) -- (Y9);
	\draw[lightgray,->] (Z10) -- (Y10);
	\draw[lightgray,->] (Z11) -- (Y11);
	%\draw[lightgray,->] (Z12) -- (Y12);

	\draw[lightgray,->] (Y1) -- (Z2);
	\draw[lightgray,->] (Y2) -- (Z3);
	\draw[lightgray,->] (Y3) -- (Z4);
	\draw[lightgray,->] (Y4) -- (Z5);
	\draw[->] (Y5) -- (Z6);
	\draw[->] (Y6) -- (Z7);
	\draw[lightgray,->] (Y7) -- (Z8);
	\draw[lightgray,->] (Y8) -- (Z9);
	\draw[lightgray,->] (Y9) -- (Z10);
	\draw[lightgray,->] (Y10) -- (Z11);
	\draw[lightgray,->] (Y11) -- (Z12);
	
	%arrows 2nd to 3rd
	\draw[lightgray,<-] (Y1) -- (S1);
	\draw[lightgray,<-] (Y2) -- (S2);
	\draw[lightgray,<-] (Y3) -- (S3);
	\draw[lightgray,<-] (Y4) -- (S4);
	\draw[<-] (Y5) -- (S5);
	\draw[<-] (Y6) -- (S6);
	\draw[<-] (Y7) -- (S7);
	\draw[lightgray,<-] (Y8) -- (S8);
	\draw[lightgray,<-] (Y9) -- (S9);
	\draw[lightgray,<-] (Y10) -- (S10);
	\draw[lightgray,<-] (Y11) -- (S11);
	%\draw[lightgray,<-] (Y12) -- (S12);
	
	\draw[lightgray,->] (Y1) -- (S2);
	\draw[lightgray,->] (Y2) -- (S3);
	\draw[lightgray,->] (Y3) -- (S4);
	\draw[lightgray,->] (Y4) -- (S5);
	\draw[->] (Y5) -- (S6);
	\draw[->] (Y6) -- (S7);
	\draw[->] (Y7) -- (S8);
	\draw[lightgray,->] (Y8) -- (S9);
	\draw[lightgray,->] (Y9) -- (S10);
	\draw[lightgray,->] (Y10) -- (S11);
	\draw[lightgray,->] (Y11) -- (S12);	
	
	%arrows 3rd to 4th
	\draw[lightgray,->] (S1) -- (B1);
	\draw[lightgray,->] (S2) -- (B2);
	\draw[lightgray,->] (S3) -- (B3);
	\draw[lightgray,->] (S4) -- (B4);
	\draw[->] (S5) -- (B5);
	\draw[->] (S6) -- (B6);
	\draw[->] (S7) -- (B7);
	\draw[->] (S8) -- (B8);
	\draw[lightgray,->] (S9) -- (B9);
	\draw[lightgray,->] (S10) -- (B10);
	\draw[lightgray,->] (S11) -- (B11);
	%\draw[lightgray,->] (S12) -- (B12);

	\draw[lightgray,->] (B1) -- (S2);
	\draw[lightgray,->] (B2) -- (S3);
	\draw[lightgray,->] (B3) -- (S4);
	\draw[->] (B4) -- (S5);
	\draw[->] (B5) -- (S6);
	\draw[->] (B6) -- (S7);
	\draw[->] (B7) -- (S8);
	\draw[lightgray,->] (B8) -- (S9);
	\draw[lightgray,->] (B9) -- (S10);
	\draw[lightgray,->] (B10) -- (S11);
	\draw[lightgray,->] (B11) -- (S12);
	
	%arrows 4th to 5th
	\draw[lightgray,<-] (B1) -- (C1);
	\draw[lightgray,<-] (B2) -- (C2);
	\draw[lightgray,<-] (B3) -- (C3);
	\draw[<-] (B4) -- (C4);
	\draw[<-] (B5) -- (C5);
	\draw[<-] (B6) -- (C6);
	\draw[<-] (B7) -- (C7);
	\draw[<-] (B8) -- (C8);
	\draw[lightgray,<-] (B9) -- (C9);
	\draw[lightgray,<-] (B10) -- (C10);
	\draw[lightgray,<-] (B11) -- (C11);
	%\draw[lightgray,<-] (B12) -- (C12);
	
	\draw[lightgray,->] (B1) -- (C2);
	\draw[lightgray,->] (B2) -- (C3);
	\draw[lightgray,->] (B3) -- (C4);
	\draw[->] (B4) -- (C5);
	\draw[->] (B5) -- (C6);
	\draw[->] (B6) -- (C7);
	\draw[->] (B7) -- (C8);
	\draw[->] (B8) -- (C9);
	\draw[lightgray,->] (B9) -- (C10);
	\draw[lightgray,->] (B10) -- (C11);
	\draw[lightgray,->] (B11) -- (C12);	
	
	%arrows 5th to 6th
	\draw[lightgray,->] (C1) -- (D1);
	\draw[lightgray,->] (C2) -- (D2);
	\draw[lightgray,->] (C3) -- (D3);
	\draw[->] (C4) -- (D4);
	\draw[->] (C5) -- (D5);
	\draw[->] (C6) -- (D6);
	\draw[->] (C7) -- (D7);
	\draw[->] (C8) -- (D8);
	\draw[->] (C9) -- (D9);
	\draw[lightgray,->] (C10) -- (D10);
	\draw[lightgray,->] (C11) -- (D11);
	%\draw[lightgray,->] (C12) -- (D12);
	
	\draw[lightgray,->] (D1) -- (C2);
	\draw[lightgray,->] (D2) -- (C3);
	\draw[->] (D3) -- (C4);
	\draw[->] (D4) -- (C5);
	\draw[->] (D5) -- (C6);
	\draw[->] (D6) -- (C7);
	\draw[->] (D7) -- (C8);
	\draw[->] (D8) -- (C9);
	\draw[lightgray,->] (D9) -- (C10);
	\draw[lightgray,->] (D10) -- (C11);
	\draw[lightgray,->] (D11) -- (C12);
	
	%arrows 6th to 7th
	\draw[lightgray,<-] (D1) -- (E1);
	\draw[lightgray,<-] (D2) -- (E2);
	\draw[<-] (D3) -- (E3);
	\draw[<-] (D4) -- (E4);
	\draw[<-] (D5) -- (E5);
	\draw[<-] (D6) -- (E6);
	\draw[<-] (D7) -- (E7);
	\draw[<-] (D8) -- (E8);
	\draw[<-] (D9) -- (E9);
	\draw[lightgray,<-] (D10) -- (E10);
	\draw[lightgray,<-] (D11) -- (E11);
	%\draw[lightgray,<-] (D12) -- (E12);
	
	\draw[lightgray,->] (D1) -- (E2);
	\draw[lightgray,->] (D2) -- (E3);
	\draw[->] (D3) -- (E4);
	\draw[->] (D4) -- (E5);
	\draw[->] (D5) -- (E6);
	\draw[->] (D6) -- (E7);
	\draw[->] (D7) -- (E8);
	\draw[->] (D8) -- (E9);
	\draw[->] (D9) -- (E10);
	\draw[lightgray,->] (D10) -- (E11);
	\draw[lightgray,->] (D11) -- (E12);
	
	%arrows 7th to 8th
	\draw[lightgray,->] (E1) -- (F1);
	\draw[lightgray,->] (E2) -- (F2);
	\draw[->] (E3) -- (F3);
	\draw[->] (E4) -- (F4);
	\draw[->] (E5) -- (F5);
	\draw[->] (E6) -- (F6);
	\draw[->] (E7) -- (F7);
	\draw[->] (E8) -- (F8);
	\draw[->] (E9) -- (F9);
	\draw[->] (E10) -- (F10);
	\draw[lightgray,->] (E11) -- (F11);
	%\draw[lightgray,->] (E12) -- (F12);
	
	\draw[lightgray,->] (F1) -- (E2);
	\draw[->] (F2) -- (E3);
	\draw[->] (F3) -- (E4);
	\draw[->] (F4) -- (E5);
	\draw[->] (F5) -- (E6);
	\draw[->] (F6) -- (E7);
	\draw[->] (F7) -- (E8);
	\draw[->] (F8) -- (E9);
	\draw[->] (F9) -- (E10);
	\draw[lightgray,->] (F10) -- (E11);
	\draw[lightgray,->] (F11) -- (E12);

	%arrows 8th to 9th
	\draw[lightgray,<-] (F1) -- (G1);
	\draw[<-] (F2) -- (G2);
	\draw[<-] (F3) -- (G3);
	\draw[<-] (F4) -- (G4);
	\draw[<-] (F5) -- (G5);
	\draw[<-] (F6) -- (G6);
	\draw[<-] (F7) -- (G7);
	\draw[<-] (F8) -- (G8);
	\draw[<-] (F9) -- (G9);
	\draw[<-] (F10) -- (G10);
	\draw[lightgray,<-] (F11) -- (G11);
	%\draw[lightgray,<-] (F12) -- (G12);
	
	\draw[lightgray,->] (F1) -- (G2);
	\draw[->] (F2) -- (G3);
	\draw[->] (F3) -- (G4);
	\draw[->] (F4) -- (G5);
	\draw[->] (F5) -- (G6);
	\draw[->] (F6) -- (G7);
	\draw[->] (F7) -- (G8);
	\draw[->] (F8) -- (G9);
	\draw[->] (F9) -- (G10);
	\draw[->] (F10) -- (G11);
	\draw[lightgray,->] (F11) -- (G12);
	
	%arrows 9th to 10th
	\draw[lightgray,->] (G1) -- (H1);
	\draw[->] (G2) -- (H2);
	\draw[->] (G3) -- (H3);
	\draw[->] (G4) -- (H4);
	\draw[->] (G5) -- (H5);
	\draw[->] (G6) -- (H6);
	\draw[->] (G7) -- (H7);
	\draw[->] (G8) -- (H8);
	\draw[->] (G9) -- (H9);
	\draw[->] (G10) -- (H10);
	\draw[->] (G11) -- (H11);
	%\draw[lightgray,->] (G12) -- (H12);
	
	\draw[->] (H1) -- (G2);
	\draw[->] (H2) -- (G3);
	\draw[->] (H3) -- (G4);
	\draw[->] (H4) -- (G5);
	\draw[->] (H5) -- (G6);
	\draw[->] (H6) -- (G7);
	\draw[->] (H7) -- (G8);
	\draw[->] (H8) -- (G9);
	\draw[->] (H9) -- (G10);
	\draw[->] (H10) -- (G11);
	\draw[lightgray,->] (H11) -- (G12);
	
	%arrows 10th to 11th
	\draw[<-] (H1) -- (I1);
	\draw[<-] (H2) -- (I2);
	\draw[<-] (H3) -- (I3);
	\draw[<-] (H4) -- (I4);
	\draw[<-] (H5) -- (I5);
	\draw[<-] (H6) -- (I6);
	\draw[<-] (H7) -- (I7);
	\draw[<-] (H8) -- (I8);
	\draw[<-] (H9) -- (I9);
	\draw[<-] (H10) -- (I10);
	\draw[<-] (H11) -- (I11);
	%\draw[lightgray,<-] (F12) -- (I12);
	
	\draw[->] (H1) -- (I2);
	\draw[->] (H2) -- (I3);
	\draw[->] (H3) -- (I4);
	\draw[->] (H4) -- (I5);
	\draw[->] (H5) -- (I6);
	\draw[->] (H6) -- (I7);
	\draw[->] (H7) -- (I8);
	\draw[->] (H8) -- (I9);
	\draw[->] (H9) -- (I10);
	\draw[->] (H10) -- (I11);
	\draw[->] (H11) -- (I12);

\end{tikzpicture}
\end{center}

\noindent Consider the cluster-tilted algebra $B_T = \End_{\CC_{A_{11}}}(T)$. 
Then $B_T=kQ \big/I$, where $Q$ is the quiver
\vskip-2ex
\[\xymatrix@!C=0pt@!R=0pt{
&&8\ar[rd]&&&&9\ar[rd]&\\
Q: &7\ar[ru]&&3\ar[ll] \ar[rd]&&2\ar[ll] \ar[ru]&&10 \ar[ll] \ar[d]  \\
&&&&1\ar[ru]\ar[ld]&&&11 \\
&6&&5\ar[ll] \ar[rr]&&4 \ar[lu]}\]
and $I$ is the ideal generated by the directed paths of length 2 which are part of the same 3-cycle. 
We refer the reader to \cite{BuanMarshReiten} for a detailed description of cluster-tilted algebras of Dynkin type $A$.

We can view $\modcat(B_T)$ as a subcategory of $\CC_{A_{11}}$ and label the indecomposable objects in 
$\CC_{A_{11}}$ by modules and shifts of projective modules respectively:

\begin{center}
\begin{tikzpicture}[scale=0.65, transform shape]

	%first row
	\node[lightgray] (Z1) at (0,10) {$\begin{smallmatrix} 8\\3\\1\\5\\6 \end{smallmatrix}$};
	\node[lightgray] (Z2) at (2,10) {$P_6[1]$};
	\node[lightgray] (Z3) at (4,10) {$6$};
	\node[lightgray] (Z4) at (6,10) {$\begin{smallmatrix}5\\4 \end{smallmatrix}$};
	\node[lightgray] (Z5) at (8,10) {$P_4[1]$};
	\node[] (Z6) at (10,10) {$\begin{smallmatrix} 4\\1\\2\\9 \end{smallmatrix}$};
	\node (Z7) at (12,10) {$P_9[1]$};
	\node[lightgray] (Z8) at (14,10) {$\begin{smallmatrix} 9\\10\\11 \end{smallmatrix}$};
	\node[lightgray] (Z9) at (16,10) {$P_{11}[1]$};
	\node[lightgray, ] (Z10) at (18,10) {$11$};
	\node[lightgray] (Z11) at (20,10) {$\begin{smallmatrix}10\\2\\3\\7\end{smallmatrix}$};
	\node[lightgray] (Z12) at (22,10) {$P_7[1]$};
	
	%2nd row
	\node[lightgray] (Y1) at (1,9) {$\begin{smallmatrix} 8\\3\\1\\5 \end{smallmatrix}$};
	\node[lightgray] (Y2) at (3,9) {$P_5[1]$};
	\node[lightgray, ] (Y3) at (5,9) {$\begin{smallmatrix} &5&\\6&&4 \end{smallmatrix}$};
	\node[lightgray] (Y4) at (7,9) {$5$};
	\node (Y5) at (9,9) {$\begin{smallmatrix} 1\\2\\9 \end{smallmatrix}$};
	\node[] (Y6) at (11,9) {$\begin{smallmatrix}4\\1\\2 \end{smallmatrix}$};
	\node (Y7) at (13,9) {$\begin{smallmatrix} 10\\11 \end{smallmatrix}$};
	\node[lightgray] (Y8) at (15,9) {$\begin{smallmatrix}9\\10\end{smallmatrix}$};
	\node[lightgray] (Y9) at (17,9) {$P_{10}[1]$}; 
	\node[lightgray, ] (Y10) at (19,9) {$\begin{smallmatrix}&10\\2&11\\3&\\7& \end{smallmatrix}$};
	\node[lightgray] (Y11) at (21,9) {$\begin{smallmatrix} 10\\2\\3 \end{smallmatrix}$};
	%\node[lightgray] (Y12) at (23,9) {$\begin{smallmatrix} \bullet \end{smallmatrix}$};
	
	%3rd row
	\node[lightgray] (S1) at (0,8) {$\begin{smallmatrix} 3\\1\\5 \end{smallmatrix}$};
	\node[lightgray] (S2) at (2,8) {$\begin{smallmatrix} 8\\3\\1 \end{smallmatrix}$};
	\node[lightgray] (S3) at (4,8) {$4$};
	\node[lightgray, ] (S4) at (6,8) {$\begin{smallmatrix} 5\\6 \end{smallmatrix}$};
	\node[] (S5) at (8,8) {$\begin{smallmatrix}&1&\\2&&5\\9&& \end{smallmatrix}$};
	\node (S6) at (10,8) {$\begin{smallmatrix} 1\\2 \end{smallmatrix}$};
	\node (S7) at (12,8) {$\begin{smallmatrix} &&&4\\&10&&1\\11&&2& \end{smallmatrix}$};
	\node[] (S8) at (14,8) {$10$};
	\node[lightgray, ] (S9) at (16,8) {$9$};
	\node[lightgray] (S10) at (18,8) {$\begin{smallmatrix} 2\\3\\7 \end{smallmatrix}$};
	\node[lightgray] (S11) at (20,8) {$\begin{smallmatrix} &10\\2&11\\3& \end{smallmatrix}$};
	\node[lightgray] (S12) at (22,8) {$\begin{smallmatrix} &&&10\\8&&2&\\&3&& \end{smallmatrix}$};
	
	%4th row
	\node[lightgray] (B1) at (1,7) {$\begin{smallmatrix} 3\\1 \end{smallmatrix}$};
	\node[lightgray] (B2) at (3,7) {$\begin{smallmatrix} &&&8\\4&&3&\\&1&& \end{smallmatrix}$};
	\node[lightgray] (B3) at (5,7) {$P_1[1]$};
	\node[] (B4) at (7,7) {$\begin{smallmatrix} &1&\\2&&5\\9&&6 \end{smallmatrix}$};
	\node (B5) at (9,7) {$\begin{smallmatrix} &1&\\2&&5 \end{smallmatrix}$};
	\node (B6) at (11,7) {$\begin{smallmatrix} &10&&1\\11&&2 \end{smallmatrix}$};
	\node (B7) at (13,7) {$\begin{smallmatrix} &&4\\10&&1\\&2& \end{smallmatrix}$};
	\node (B8) at (15,7) {$P_2[1]$};
	\node[lightgray, ] (B9) at (17,7) {$\begin{smallmatrix} &2&\\9&&3\\&&7 \end{smallmatrix}$}; 
	\node[lightgray] (B10) at (19,7) {$\begin{smallmatrix} 2\\3 \end{smallmatrix}$};
	\node[lightgray] (B11) at (21,7) {$\begin{smallmatrix} &&&10\\8&&2&11\\&3&& \end{smallmatrix}$};
	%\node[lightgray] (B12) at (23,7) {$\begin{smallmatrix} \bullet \end{smallmatrix}$};
	
	%5th row
	\node[lightgray] (C1) at (0,6) {$\begin{smallmatrix} &3&\\1&&7 \end{smallmatrix}$};
	\node[lightgray] (C2) at (2,6) {$\begin{smallmatrix} 4&&3\\&1& \end{smallmatrix}$};
	\node[lightgray] (C3) at (4,6) {$\begin{smallmatrix} 8\\3 \end{smallmatrix}$};
	\node (C4) at (6,6) {$\begin{smallmatrix} 2\\9 \end{smallmatrix}$};
	\node (C5) at (8,6) {$\begin{smallmatrix} &1&\\2&&5 \\&&6\end{smallmatrix}$};
	\node (C6) at (10,6) {$\begin{smallmatrix}&10&&1&\\11&&2&&5\end{smallmatrix}$};
	\node (C7) at (12,6) {$\begin{smallmatrix} 10&&1\\&2& \end{smallmatrix}$};
	\node (C8) at (14,6) {$\begin{smallmatrix} 4\\1 \end{smallmatrix}$};
	\node[] (C9) at (16,6) {$\begin{smallmatrix} 3\\7 \end{smallmatrix}$};
	\node[lightgray] (C10) at (18,6) {$\begin{smallmatrix} &2&\\9&&3 \end{smallmatrix}$};
	\node[lightgray] (C11) at (20,6) {$\begin{smallmatrix} 8&&2\\&3& \end{smallmatrix}$};
	\node[lightgray] (C12) at (22,6) {$\begin{smallmatrix}&10&\\2&&11 \end{smallmatrix}$};
	
	%6th row
	\node[lightgray] (D1) at (1,5) {$\begin{smallmatrix} 4&&3&\\&1&&7 \end{smallmatrix}$};
	\node[lightgray] (D2) at (3,5) {$3$};
	\node (D3) at (5,5) {$\begin{smallmatrix} &2&8\\9&3&\end{smallmatrix}$};
	\node[] (D4) at (7,5) {$2$};
	\node (D5) at (9,5) {$\begin{smallmatrix} 10&&1&\\11&2&&5&\\&&&6 \end{smallmatrix}$};
	\node (D6) at (11,5) {$\begin{smallmatrix}10&&1&\\&2&&5 \end{smallmatrix}$};
	\node (D7) at (13,5) {$1$};
	\node (D8) at (15,5) {$\begin{smallmatrix} 4&&3&\\&1&&7\end{smallmatrix}$};
	\node (D9) at (17,5) {$3$}; 
	\node[lightgray] (D10) at (19,5) {$\begin{smallmatrix} &2&&8\\9&&3& \end{smallmatrix}$};
	\node[lightgray] (D11) at (21,5) {$2$};
	%\node[lightgray] (D12) at (23,5) {$\begin{smallmatrix} \bullet \end{smallmatrix}$};
	
	%7th row
	\node[lightgray] (E1) at (0,4) {$\begin{smallmatrix} 4\\1 \end{smallmatrix}$};
	\node[lightgray] (E2) at (2,4) {$\begin{smallmatrix} 3\\7 \end{smallmatrix}$};
	\node (E3) at (4,4) {$\begin{smallmatrix} &2&\\9&&3 \end{smallmatrix}$};
	\node (E4) at (6,4) {$\begin{smallmatrix} 8&&2\\&3& \end{smallmatrix}$};
	\node (E5) at (8,4) {$\begin{smallmatrix} &10\\2&11 \end{smallmatrix}$};
	\node (E6) at (10,4) {$\begin{smallmatrix} 10&&1&\\&2&&5\\&&&6 \end{smallmatrix}$};
	\node (E7) at (12,4) {$\begin{smallmatrix} 1\\5 \end{smallmatrix}$};
	\node (E8) at (14,4) {$\begin{smallmatrix} &3&\\1&&7 \end{smallmatrix}$};
	\node (E9) at (16,4) {$\begin{smallmatrix} 4&&3\\&1& \end{smallmatrix}$};
	\node (E10) at (18,4) {$\begin{smallmatrix} 8\\3 \end{smallmatrix}$};
	\node[lightgray] (E11) at (20,4) {$\begin{smallmatrix} 2\\9 \end{smallmatrix}$};
	\node[lightgray] (E12) at (22,4) {$\begin{smallmatrix} &1&\\2&&5 \\&&6\end{smallmatrix}$};
	
	%8th row
	\node[lightgray] (F1) at (1,3) {$P_2[1]$};
	\node[] (F2) at (3,3) {$\begin{smallmatrix} &2&\\9&&3\\&&7 \end{smallmatrix}$};
	\node[] (F3) at (5,3) {$\begin{smallmatrix} 2\\3 \end{smallmatrix}$};
	\node (F4) at (7,3) {$\begin{smallmatrix} &&&10\\8&&2&11\\&3&& \end{smallmatrix}$};
	\node (F5) at (9,3) {$\begin{smallmatrix}10\\2 \end{smallmatrix}$};
	\node (F6) at (11,3) {$\begin{smallmatrix} 1\\5\\6 \end{smallmatrix}$};
	\node (F7) at (13,3) {$\begin{smallmatrix} 3&\\1&7\\5&\end{smallmatrix}$};
	\node (F8) at (15,3) {$\begin{smallmatrix} 3\\1\end{smallmatrix}$};
	\node[] (F9) at (17,3) {$\begin{smallmatrix} &&8\\4&&3\\&1& \end{smallmatrix}$};
	\node (F10) at (19,3) {$P_1[1]$};
	\node[lightgray, ] (F11) at (21,3) {$\begin{smallmatrix} &1&\\2&&5\\9&&6 \end{smallmatrix}$};
	%\node[lightgray, ] (F12) at (23,3) {$\begin{smallmatrix} \bullet \end{smallmatrix}$}; 
	
	%9th row
	\node[lightgray] (G1) at (0,2) {$10$};
	\node (G2) at (2,2) {$9$};
	\node (G3) at (4,2) {$\begin{smallmatrix} 2\\3\\7 \end{smallmatrix}$};
	\node (G4) at (6,2) {$\begin{smallmatrix}&10\\2&11\\3& \end{smallmatrix}$};
	\node (G5) at (8,2) {$\begin{smallmatrix} &&10\\8&&2\\&3& \end{smallmatrix}$};
	\node[] (G6) at (10,2) {$P_3[1]$};
	\node[] (G7) at (12,2) {$\begin{smallmatrix} &3&\\1&&7 \\5&&\\6&& \end{smallmatrix}$};
	\node (G8) at (14,2) {$\begin{smallmatrix} 3\\1\\5 \end{smallmatrix}$};
	\node (G9) at (16,2) {$\begin{smallmatrix} 8\\3\\1 \end{smallmatrix}$};
	\node (G10) at (18,2) {$4$};
	\node (G11) at (20,2) {$\begin{smallmatrix} 5\\6 \end{smallmatrix}$};
	\node[lightgray] (G12) at (22,2) {$\begin{smallmatrix} &1&\\2&&5\\9&& \end{smallmatrix}$};
	
	%10th row
	\node (H1) at (1,1) {$\begin{smallmatrix} 9\\10 \end{smallmatrix}$};
	\node[] (H2) at (3,1) {$P_{10}[1]$};
	\node[] (H3) at (5,1) {$\begin{smallmatrix} &10\\2&11\\3&\\7& \end{smallmatrix}$};
	\node[] (H4) at (7,1) {$\begin{smallmatrix} 10\\2\\3 \end{smallmatrix}$};
	\node (H5) at (9,1) {$8$};
	\node[] (H6) at (11,1) {$7$};
	\node (H7) at (13,1) {$\begin{smallmatrix} 3\\1\\5\\6 \end{smallmatrix}$};
	\node[] (H8) at (15,1) {$\begin{smallmatrix} 8\\3\\1\\5 \end{smallmatrix}$};
	\node (H9) at (17,1) {$P_5[1]$}; 
	\node[] (H10) at (19,1) {$\begin{smallmatrix} &5&\\6&&4 \end{smallmatrix}$};
	\node[] (H11) at (21,1) {$5$};
	%\node[lightgray] (H12) at (23,1) {$\begin{smallmatrix} \bullet \end{smallmatrix}$};
	
	%11th row
	\node[] (I1) at (0,0) {$\begin{smallmatrix} 9\\10\\11 \end{smallmatrix}$};
	\node[] (I2) at (2,0) {$P_{11}[1]$};
	\node[] (I3) at (4,0) {$11$};
	\node[] (I4) at (6,0) {$\begin{smallmatrix} 10\\2\\3\\7 \end{smallmatrix}$};
	\node[] (I5) at (8,0) {$P_7[1]$};
	\node[] (I6) at (10,0) {$\begin{smallmatrix} 7\\8 \end{smallmatrix}$};
	\node[] (I7) at (12,0) {$P_8[1]$};
	\node[] (I8) at (14,0) {$\begin{smallmatrix} 8\\3\\1\\5\\6 \end{smallmatrix}$};
	\node[]  (I9) at (16,0) {$P_6[1]$}; 
	\node[] (I10) at (18,0) {$6$};
	\node (I11) at (20,0) {$\begin{smallmatrix} 5\\4 \end{smallmatrix}$};
	\node[] (I12) at (22,0) {$P_4[1]$};
	
	%arrows 1st to 2nd
	\draw[lightgray, ->] (Z1) -- (Y1);
	\draw[lightgray,->] (Z2) -- (Y2);
	\draw[lightgray,->] (Z3) -- (Y3);
	\draw[lightgray,->] (Z4) -- (Y4);
	\draw[lightgray,->] (Z5) -- (Y5);
	\draw[->] (Z6) -- (Y6);
	\draw[->] (Z7) -- (Y7);
	\draw[lightgray,->] (Z8) -- (Y8);
	\draw[lightgray,->] (Z9) -- (Y9);
	\draw[lightgray,->] (Z10) -- (Y10);
	\draw[lightgray,->] (Z11) -- (Y11);
	%\draw[lightgray,->] (Z12) -- (Y12);

	\draw[lightgray,->] (Y1) -- (Z2);
	\draw[lightgray,->] (Y2) -- (Z3);
	\draw[lightgray,->] (Y3) -- (Z4);
	\draw[lightgray,->] (Y4) -- (Z5);
	\draw[->] (Y5) -- (Z6);
	\draw[->] (Y6) -- (Z7);
	\draw[lightgray,->] (Y7) -- (Z8);
	\draw[lightgray,->] (Y8) -- (Z9);
	\draw[lightgray,->] (Y9) -- (Z10);
	\draw[lightgray,->] (Y10) -- (Z11);
	\draw[lightgray,->] (Y11) -- (Z12);
	
	%arrows 2nd to 3rd
	\draw[lightgray,<-] (Y1) -- (S1);
	\draw[lightgray,<-] (Y2) -- (S2);
	\draw[lightgray,<-] (Y3) -- (S3);
	\draw[lightgray,<-] (Y4) -- (S4);
	\draw[<-] (Y5) -- (S5);
	\draw[<-] (Y6) -- (S6);
	\draw[<-] (Y7) -- (S7);
	\draw[lightgray,<-] (Y8) -- (S8);
	\draw[lightgray,<-] (Y9) -- (S9);
	\draw[lightgray,<-] (Y10) -- (S10);
	\draw[lightgray,<-] (Y11) -- (S11);
	%\draw[lightgray,<-] (Y12) -- (S12);
	
	\draw[lightgray,->] (Y1) -- (S2);
	\draw[lightgray,->] (Y2) -- (S3);
	\draw[lightgray,->] (Y3) -- (S4);
	\draw[lightgray,->] (Y4) -- (S5);
	\draw[->] (Y5) -- (S6);
	\draw[->] (Y6) -- (S7);
	\draw[->] (Y7) -- (S8);
	\draw[lightgray,->] (Y8) -- (S9);
	\draw[lightgray,->] (Y9) -- (S10);
	\draw[lightgray,->] (Y10) -- (S11);
	\draw[lightgray,->] (Y11) -- (S12);	
	
	%arrows 3rd to 4th
	\draw[lightgray,->] (S1) -- (B1);
	\draw[lightgray,->] (S2) -- (B2);
	\draw[lightgray,->] (S3) -- (B3);
	\draw[lightgray,->] (S4) -- (B4);
	\draw[->] (S5) -- (B5);
	\draw[->] (S6) -- (B6);
	\draw[->] (S7) -- (B7);
	\draw[->] (S8) -- (B8);
	\draw[lightgray,->] (S9) -- (B9);
	\draw[lightgray,->] (S10) -- (B10);
	\draw[lightgray,->] (S11) -- (B11);
	%\draw[lightgray,->] (S12) -- (B12);

	\draw[lightgray,->] (B1) -- (S2);
	\draw[lightgray,->] (B2) -- (S3);
	\draw[lightgray,->] (B3) -- (S4);
	\draw[->] (B4) -- (S5);
	\draw[->] (B5) -- (S6);
	\draw[->] (B6) -- (S7);
	\draw[->] (B7) -- (S8);
	\draw[lightgray,->] (B8) -- (S9);
	\draw[lightgray,->] (B9) -- (S10);
	\draw[lightgray,->] (B10) -- (S11);
	\draw[lightgray,->] (B11) -- (S12);
	
	%arrows 4th to 5th
	\draw[lightgray,<-] (B1) -- (C1);
	\draw[lightgray,<-] (B2) -- (C2);
	\draw[lightgray,<-] (B3) -- (C3);
	\draw[<-] (B4) -- (C4);
	\draw[<-] (B5) -- (C5);
	\draw[<-] (B6) -- (C6);
	\draw[<-] (B7) -- (C7);
	\draw[<-] (B8) -- (C8);
	\draw[lightgray,<-] (B9) -- (C9);
	\draw[lightgray,<-] (B10) -- (C10);
	\draw[lightgray,<-] (B11) -- (C11);
	%\draw[lightgray,<-] (B12) -- (C12);
	
	\draw[lightgray,->] (B1) -- (C2);
	\draw[lightgray,->] (B2) -- (C3);
	\draw[lightgray,->] (B3) -- (C4);
	\draw[->] (B4) -- (C5);
	\draw[->] (B5) -- (C6);
	\draw[->] (B6) -- (C7);
	\draw[->] (B7) -- (C8);
	\draw[->] (B8) -- (C9);
	\draw[lightgray,->] (B9) -- (C10);
	\draw[lightgray,->] (B10) -- (C11);
	\draw[lightgray,->] (B11) -- (C12);	
	
	%arrows 5th to 6th
	\draw[lightgray,->] (C1) -- (D1);
	\draw[lightgray,->] (C2) -- (D2);
	\draw[lightgray,->] (C3) -- (D3);
	\draw[->] (C4) -- (D4);
	\draw[->] (C5) -- (D5);
	\draw[->] (C6) -- (D6);
	\draw[->] (C7) -- (D7);
	\draw[->] (C8) -- (D8);
	\draw[->] (C9) -- (D9);
	\draw[lightgray,->] (C10) -- (D10);
	\draw[lightgray,->] (C11) -- (D11);
	%\draw[lightgray,->] (C12) -- (D12);
	
	\draw[lightgray,->] (D1) -- (C2);
	\draw[lightgray,->] (D2) -- (C3);
	\draw[->] (D3) -- (C4);
	\draw[->] (D4) -- (C5);
	\draw[->] (D5) -- (C6);
	\draw[->] (D6) -- (C7);
	\draw[->] (D7) -- (C8);
	\draw[->] (D8) -- (C9);
	\draw[lightgray,->] (D9) -- (C10);
	\draw[lightgray,->] (D10) -- (C11);
	\draw[lightgray,->] (D11) -- (C12);
	
	%arrows 6th to 7th
	\draw[lightgray,<-] (D1) -- (E1);
	\draw[lightgray,<-] (D2) -- (E2);
	\draw[<-] (D3) -- (E3);
	\draw[<-] (D4) -- (E4);
	\draw[<-] (D5) -- (E5);
	\draw[<-] (D6) -- (E6);
	\draw[<-] (D7) -- (E7);
	\draw[<-] (D8) -- (E8);
	\draw[<-] (D9) -- (E9);
	\draw[lightgray,<-] (D10) -- (E10);
	\draw[lightgray,<-] (D11) -- (E11);
	%\draw[lightgray,<-] (D12) -- (E12);
	
	\draw[lightgray,->] (D1) -- (E2);
	\draw[lightgray,->] (D2) -- (E3);
	\draw[->] (D3) -- (E4);
	\draw[->] (D4) -- (E5);
	\draw[->] (D5) -- (E6);
	\draw[->] (D6) -- (E7);
	\draw[->] (D7) -- (E8);
	\draw[->] (D8) -- (E9);
	\draw[->] (D9) -- (E10);
	\draw[lightgray,->] (D10) -- (E11);
	\draw[lightgray,->] (D11) -- (E12);
	
	%arrows 7th to 8th
	\draw[lightgray,->] (E1) -- (F1);
	\draw[lightgray,->] (E2) -- (F2);
	\draw[->] (E3) -- (F3);
	\draw[->] (E4) -- (F4);
	\draw[->] (E5) -- (F5);
	\draw[->] (E6) -- (F6);
	\draw[->] (E7) -- (F7);
	\draw[->] (E8) -- (F8);
	\draw[->] (E9) -- (F9);
	\draw[->] (E10) -- (F10);
	\draw[lightgray,->] (E11) -- (F11);
	%\draw[lightgray,->] (E12) -- (F12);
	
	\draw[lightgray,->] (F1) -- (E2);
	\draw[->] (F2) -- (E3);
	\draw[->] (F3) -- (E4);
	\draw[->] (F4) -- (E5);
	\draw[->] (F5) -- (E6);
	\draw[->] (F6) -- (E7);
	\draw[->] (F7) -- (E8);
	\draw[->] (F8) -- (E9);
	\draw[->] (F9) -- (E10);
	\draw[lightgray,->] (F10) -- (E11);
	\draw[lightgray,->] (F11) -- (E12);

	%arrows 8th to 9th
	\draw[lightgray,<-] (F1) -- (G1);
	\draw[<-] (F2) -- (G2);
	\draw[<-] (F3) -- (G3);
	\draw[<-] (F4) -- (G4);
	\draw[<-] (F5) -- (G5);
	\draw[<-] (F6) -- (G6);
	\draw[<-] (F7) -- (G7);
	\draw[<-] (F8) -- (G8);
	\draw[<-] (F9) -- (G9);
	\draw[<-] (F10) -- (G10);
	\draw[lightgray,<-] (F11) -- (G11);
	%\draw[lightgray,<-] (F12) -- (G12);
	
	\draw[lightgray,->] (F1) -- (G2);
	\draw[->] (F2) -- (G3);
	\draw[->] (F3) -- (G4);
	\draw[->] (F4) -- (G5);
	\draw[->] (F5) -- (G6);
	\draw[->] (F6) -- (G7);
	\draw[->] (F7) -- (G8);
	\draw[->] (F8) -- (G9);
	\draw[->] (F9) -- (G10);
	\draw[->] (F10) -- (G11);
	\draw[lightgray,->] (F11) -- (G12);
	
	%arrows 9th to 10th
	\draw[lightgray,->] (G1) -- (H1);
	\draw[->] (G2) -- (H2);
	\draw[->] (G3) -- (H3);
	\draw[->] (G4) -- (H4);
	\draw[->] (G5) -- (H5);
	\draw[->] (G6) -- (H6);
	\draw[->] (G7) -- (H7);
	\draw[->] (G8) -- (H8);
	\draw[->] (G9) -- (H9);
	\draw[->] (G10) -- (H10);
	\draw[->] (G11) -- (H11);
	%\draw[lightgray,->] (G12) -- (H12);
	
	\draw[->] (H1) -- (G2);
	\draw[->] (H2) -- (G3);
	\draw[->] (H3) -- (G4);
	\draw[->] (H4) -- (G5);
	\draw[->] (H5) -- (G6);
	\draw[->] (H6) -- (G7);
	\draw[->] (H7) -- (G8);
	\draw[->] (H8) -- (G9);
	\draw[->] (H9) -- (G10);
	\draw[->] (H10) -- (G11);
	\draw[lightgray,->] (H11) -- (G12);
	
	%arrows 10th to 11th
	\draw[<-] (H1) -- (I1);
	\draw[<-] (H2) -- (I2);
	\draw[<-] (H3) -- (I3);
	\draw[<-] (H4) -- (I4);
	\draw[<-] (H5) -- (I5);
	\draw[<-] (H6) -- (I6);
	\draw[<-] (H7) -- (I7);
	\draw[<-] (H8) -- (I8);
	\draw[<-] (H9) -- (I9);
	\draw[<-] (H10) -- (I10);
	\draw[<-] (H11) -- (I11);
	%\draw[lightgray,<-] (F12) -- (I12);
	
	\draw[->] (H1) -- (I2);
	\draw[->] (H2) -- (I3);
	\draw[->] (H3) -- (I4);
	\draw[->] (H4) -- (I5);
	\draw[->] (H5) -- (I6);
	\draw[->] (H6) -- (I7);
	\draw[->] (H7) -- (I8);
	\draw[->] (H8) -- (I9);
	\draw[->] (H9) -- (I10);
	\draw[->] (H10) -- (I11);
	\draw[->] (H11) -- (I12);

\end{tikzpicture}
\end{center}

\noindent The specialized CC-map
replaces each vertex labelled by a module, by the number of its submodules and the shifts of projectives by 1s. Adding in the first two and last two rows of 0s and 1s gives rise to the associated frieze $F(T)$:
\vskip-2ex
\begin{center}
\begin{tikzpicture}[scale=0.6, transform shape]

	%0s
	\node[lightgray] (1) at (0,12) {$0$};
	\node[lightgray] (2) at (2,12) {$0$};
	\node[lightgray] (3) at (4,12) {$0$};
	\node[lightgray] (4) at (6,12) {$0$};
	\node[lightgray] (5) at (8,12) {$0$};
	\node[lightgray] (6) at (10,12) {$0$};
	\node[lightgray] (7) at (12,12) {$0$};
	\node[lightgray] (8) at (14,12) {$0$};
	\node[lightgray] (9) at (16,12) {$0$};
	\node[lightgray] (10) at (18,12) {$0$};
	\node[lightgray] (11) at (20,12) {$0$};
	\node[lightgray] (12) at (22,12) {$0$};
	
	%1s
	\node[lightgray] (1) at (1,11) {$1$};
	\node[lightgray] (2) at (3,11) {$1$};
	\node[lightgray] (3) at (5,11) {$1$};
	\node[lightgray] (4) at (7,11) {$1$};
	\node[lightgray] (5) at (9,11) {$1$};
	\node[lightgray] (6) at (11,11) {$1$};
	\node[lightgray] (7) at (13,11) {$1$};
	\node[lightgray] (8) at (15,11) {$1$};
	\node[lightgray] (9) at (17,11) {$1$}; 
	\node[lightgray] (10) at (19,11) {$1$};
	\node[lightgray] (11) at (21,11) {$1$};
	\node[lightgray] (12) at (23,11) {$1$};

	%first row
	\node[lightgray] (Z1) at (0,10) {$6$};
	\node[lightgray] (Z2) at (2,10) {$1$};
	\node[lightgray] (Z3) at (4,10) {$2$};
	\node[lightgray] (Z4) at (6,10) {$3$};
	\node[lightgray] (Z5) at (8,10) {$1$};
	\node (Z6) at (10,10) {$5$};
	\node (Z7) at (12,10) {$1$};
	\node[lightgray] (Z8) at (14,10) {$4$};
	\node[lightgray] (Z9) at (16,10) {$1$};
	\node[lightgray] (Z10) at (18,10) {$2$};
	\node[lightgray] (Z11) at (20,10) {$5$};
	\node[lightgray] (Z12) at (22,10) {$1$};
	
	%2nd row
	\node[lightgray] (Y1) at (1,9) {$5$};
	\node[lightgray] (Y2) at (3,9) {$1$};
	\node[lightgray] (Y3) at (5,9) {$5$};
	\node[lightgray] (Y4) at (7,9) {$2$};
	\node (Y5) at (9,9) {$4$};
	\node (Y6) at (11,9) {$4$};
	\node (Y7) at (13,9) {$3$};
	\node[lightgray] (Y8) at (15,9) {$3$};
	\node[lightgray] (Y9) at (17,9) {$1$}; 
	\node[lightgray] (Y10) at (19,9) {$9$};
	\node[lightgray] (Y11) at (21,9) {$4$};
	\node[lightgray] (Y12) at (23,9) {$2$};
	
	%3rd row
	\node[lightgray] (S1) at (0,8) {$4$};
	\node[lightgray] (S2) at (2,8) {$4$};
	\node[lightgray] (S3) at (4,8) {$2$};
	\node[lightgray] (S4) at (6,8) {$3$};
	\node (S5) at (8,8) {$7$};
	\node (S6) at (10,8) {$3$};
	\node (S7) at (12,8) {$11$};
	\node (S8) at (14,8) {$2$};
	\node[lightgray] (S9) at (16,8) {$2$};
	\node[lightgray] (S10) at (18,8) {$4$};
	\node[lightgray] (S11) at (20,8) {$7$};
	\node[lightgray] (S12) at (22,8) {$7$};
	
	%4th row
	\node[lightgray] (B1) at (1,7) {$3$};
	\node[lightgray] (B2) at (3,7) {$7$};
	\node[lightgray] (B3) at (5,7) {$1$};
	\node (B4) at (7,7) {$10$};
	\node (B5) at (9,7) {$5$};
	\node (B6) at (11,7) {$8$};
	\node (B7) at (13,7) {$7$};
	\node (B8) at (15,7) {$1$};
	\node[lightgray] (B9) at (17,7) {$7$}; 
	\node[lightgray] (B10) at (19,7) {$3$};
	\node[lightgray] (B11) at (21,7) {$12$};
	\node[lightgray] (B12) at (23,7) {$3$};
	%5th row
	\node[lightgray] (C1) at (0,6) {$5$};
	\node[lightgray] (C2) at (2,6) {$5$};
	\node[lightgray] (C3) at (4,6) {$3$};
	\node (C4) at (6,6) {$3$};
	\node (C5) at (8,6) {$7$};
	\node (C6) at (10,6) {$13$};
	\node (C7) at (12,6) {$5$};
	\node (C8) at (14,6) {$3$};
	\node (C9) at (16,6) {$3$};
	\node[lightgray] (C10) at (18,6) {$5$};
	\node[lightgray] (C11) at (20,6) {$5$};
	\node[lightgray] (C12) at (22,6) {$5$};
	%6th row
	\node[lightgray] (D1) at (1,5) {$8$};
	\node[lightgray] (D2) at (3,5) {$2$};
	\node (D3) at (5,5) {$8$};
	\node (D4) at (7,5) {$2$};
	\node (D5) at (9,5) {$18$};
	\node (D6) at (11,5) {$8$};
	\node (D7) at (13,5) {$2$};
	\node (D8) at (15,5) {$8$};
	\node (D9) at (17,5) {$2$}; 
	\node[lightgray] (D10) at (19,5) {$8$};
	\node[lightgray] (D11) at (21,5) {$2$};
	\node[lightgray] (D12) at (23,5) {$18$};
	%7th row
	\node[lightgray] (E1) at (0,4) {$3$};
	\node[lightgray] (E2) at (2,4) {$3$};
	\node (E3) at (4,4) {$5$};
	\node (E4) at (6,4) {$5$};
	\node (E5) at (8,4) {$5$};
	\node (E6) at (10,4) {$11$};
	\node (E7) at (12,4) {$3$};
	\node (E8) at (14,4) {$5$};
	\node (E9) at (16,4) {$5$};
	\node (E10) at (18,4) {$3$};
	\node[lightgray] (E11) at (20,4) {$3$};
	\node[lightgray] (E12) at (22,4) {$7$};
	%8th row
	\node[lightgray] (F1) at (1,3) {$1$};
	\node (F2) at (3,3) {$7$};
	\node (F3) at (5,3) {$3$};
	\node (F4) at (7,3) {$12$};
	\node (F5) at (9,3) {$3$};
	\node (F6) at (11,3) {$4$};
	\node (F7) at (13,3) {$7$};
	\node (F8) at (15,3) {$3$};
	\node (F9) at (17,3) {$7$};
	\node (F10) at (19,3) {$1$};
	\node[lightgray] (F11) at (21,3) {$10$};
	\node[lightgray] (F12) at (23,3) {$5$}; 
	%9th row
	\node[lightgray] (G1) at (0,2) {$2$};
	\node (G2) at (2,2) {$2$};
	\node (G3) at (4,2) {$4$};
	\node (G4) at (6,2) {$7$};
	\node (G5) at (8,2) {$7$};
	\node (G6) at (10,2) {$1$};
	\node (G7) at (12,2) {$9$};
	\node (G8) at (14,2) {$4$};
	\node (G9) at (16,2) {$4$};
	\node (G10) at (18,2) {$2$};
	\node (G11) at (20,2) {$3$};
	\node[lightgray] (G12) at (22,2) {$7$};
	%10th row
	\node (H1) at (1,1) {$3$};
	\node (H2) at (3,1) {$1$};
	\node (H3) at (5,1) {$9$};
	\node (H4) at (7,1) {$4$};
	\node (H5) at (9,1) {$2$};
	\node (H6) at (11,1) {$2$};
	\node (H7) at (13,1) {$5$};
	\node (H8) at (15,1) {$5$};
	\node (H9) at (17,1) {$1$}; 
	\node (H10) at (19,1) {$5$};
	\node (H11) at (21,1) {$2$};
	\node[lightgray] (H12) at (23,1) {$4$};
	
	%11th row
	\node[] (1) at (0,0) {$4$};
	\node[] (2) at (2,0) {$1$};
	\node[] (3) at (4,0) {$2$};
	\node[] (4) at (6,0) {$5$};
	\node[] (5) at (8,0) {$1$};
	\node[] (6) at (10,0) {$3$};
	\node[] (7) at (12,0) {$1$};
	\node[] (8) at (14,0) {$6$};
	\node[] (9) at (16,0) {$1$};
	\node[] (10) at (18,0) {$2$};
	\node[] (11) at (20,0) {$3$};
	\node[] (12) at (22,0) {$1$};
	
	%1s
	\node[lightgray] (1) at (1,-1) {$1$};
	\node[lightgray] (2) at (3,-1) {$1$};
	\node[lightgray] (3) at (5,-1) {$1$};
	\node[lightgray] (4) at (7,-1) {$1$};
	\node[lightgray] (5) at (9,-1) {$1$};
	\node[lightgray] (6) at (11,-1) {$1$};
	\node[lightgray] (7) at (13,-1) {$1$};
	\node[lightgray] (8) at (15,-1) {$1$};
	\node[lightgray] (9) at (17,-1) {$1$}; 
	\node[lightgray] (10) at (19,-1) {$1$};
	\node[lightgray] (11) at (21,-1) {$1$};
	\node[lightgray] (12) at (23,-1) {$1$};
	
	%0s
	\node[lightgray] (1) at (0,-2) {$0$};
	\node[lightgray] (2) at (2,-2) {$0$};
	\node[lightgray] (3) at (4,-2) {$0$};
	\node[lightgray] (4) at (6,-2) {$0$};
	\node[lightgray] (5) at (8,-2) {$0$};
	\node[lightgray] (6) at (10,-2) {$0$};
	\node[lightgray] (7) at (12,-2) {$0$};
	\node[lightgray] (8) at (14,-2) {$0$};
	\node[lightgray] (9) at (16,-2) {$0$}; 
	\node[lightgray] (10) at (18,-2) {$0$};
	\node[lightgray] (11) at (20,-2) {$0$};
	\node[lightgray] (12) at (22,-2) {$0$};
\end{tikzpicture}
\end{center}
\vskip-2ex
 \end{ex}

%%%%%%%
\section{Description of the regions in the frieze}\label{sec:regions} 
%%%%%%

\noindent
{\bf The quiver of a triangulation.} 

Let $\mathcal T$ be a triangulation of an $(n+3)$-gon, and let the diagonals be labeled by $1,2,\dots,n$. 
We recall that the quiver $Q_{\mathcal T}$ of the triangulation $\mathcal T$ is defined as 
follows: the vertices of $Q_{\mathcal T}$ are the labels $\{1, 2,\dots, n\}$. 
There is an arrow $i\to j$ in case the diagonals share an endpoint and the diagonal $i$ 
can be rotated clockwise to diagonal $j$ 
(without passing through another diagonal incident with the common vertex). 
This is illustrated in Example~\ref{ex:quiver-triangul-flip} and Figure~\ref{fig:triangulations-1-2} 
below. 

Let $B=B_{\mathcal T}$ be the path algebra of $Q_{\mathcal{T}}$ modulo the relations 
arising from triangles in $Q_{\mathcal{T}}$: whenever $\alpha,\beta$ are two successive 
arrows in an oriented triangle in $Q_{\mathcal{T}}$, their composition is $0$. 
Let $P_x$ be the indecomposable projective $B$-module associated to 
the vertex $x$ and $S_x$ its simple top. 
Let $$T = \oplus_{x\in \mathcal T} P_x.$$ 
We considered $T$ as an object of the generalized cluster category $\CC=\CC_{B}$. 
Then $T$ is a cluster tilting object in $\CC$ and $B\cong\End_{\CC}(T).$
Hence $B$ is a cluster-tilted algebra, called the 
cluster-tilted algebra associated to the triangulation $\mathcal T$.
We can extend this to an object in the Frobenius category $\mathcal C_f$ by 
adding the $n+3$ projective-injective summands associated to the boundary segments $[12],[23]$ 
$,\dots,[n+3,1]$ 
of the polygon, 
with irreducible maps between the objects corresponding to diagonals/edges as follows: 
$[i-1,i+1]\to [i,i+1]$, $[i,i+1]\to [i,i+2]$ (\cite{JKS, BKM,DeLuo}). 
We denote the projective-injective associated to $[i,i+1]$ by $Q_{x_i}$. 
Let 
\[
{T}_f=\left( \oplus _{x\in\mathcal T} P_x\right) \oplus \left(Q_{x_1}\oplus \dots\oplus Q_{x_{n+3}}\right)
\]
This is a cluster tilting object of $\mathcal C_f$ in the sense of~\cite[Section 3]{FuKeller}. 
Given a $B$-module $M$, by abuse of notation, 
we denote the corresponding objects in $\mathcal{C}$ and $\mathcal{C}_f$ by $M$, 
that is $\text{Hom}_{\mathcal{C}}({T}, M)=M$.  In other words, an indecomposable object of 
$\mathcal{C}_f$ is either an indecomposable $B$-module or $Q_{x_i}$ for some 
$i\in\{1,\dots, n+3\}$ 
or of the form $P_x[1]$ for some $x \in \mathcal T$. 

The {\em frieze $F(\mathcal T)$ of the triangulation $\mathcal T$} 
is the frieze pattern $F(T)$ for 
$T$ the cluster tilting object associated to $\mathcal T$.

%%%%%%
\subsection{Diagonal defines quadrilateral}
%%%%%%
$\ $

\begin{figure}
\begin{center}
\begin{minipage}{65mm} \includegraphics[width=45mm]{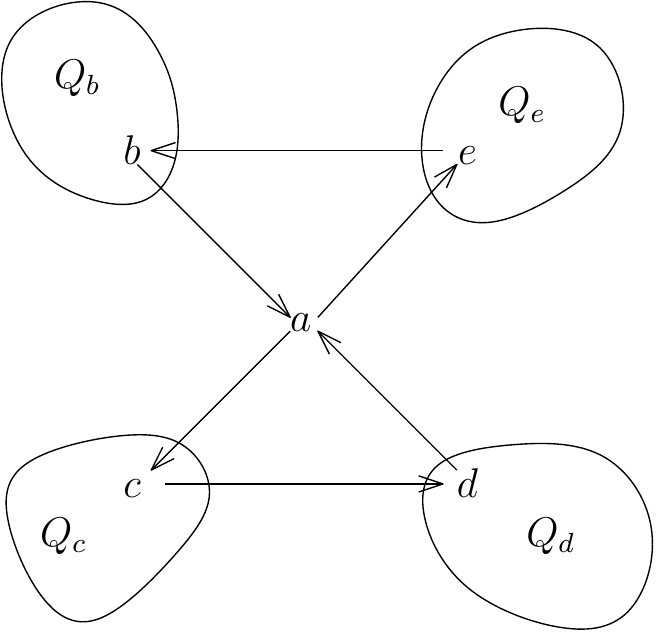}
 \caption{Regions in quiver.}\label{fig:quiver-divided}
\end{minipage}
\hfil
\begin{minipage}{65mm} \includegraphics[width=45mm]{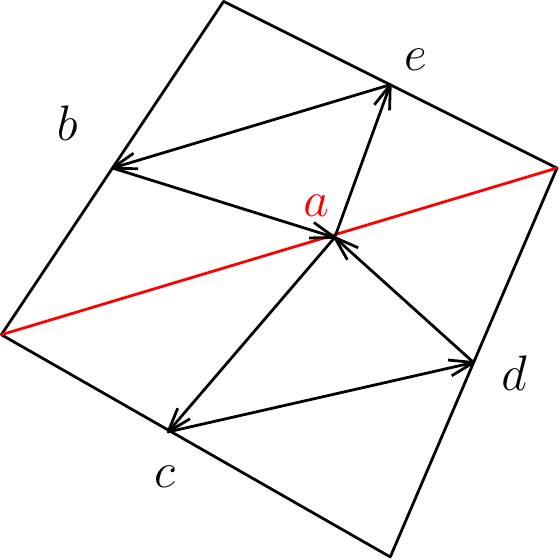} 
\caption{Triangulation around $a$.}
\label{fig:local-triang}
\end{minipage}
\end{center}
\end{figure}

Let $a$ be a diagonal in the triangulation, $a\in\{1,2,\dots, n\}$. 
This diagonal uniquely defines a quadrilateral formed by diagonals or 
boundary segments. Label them $b,c,d,e$ as in Figure~\ref{fig:local-triang}.  

%%%%%%
\subsection{Diagonal defines two rays} \label{ssec:diagonal-rays}
%%%%%%
$\ $

Consider the entry $1$ of the frieze corresponding to $a$. There are two rays passing through it. We go along 
these rays forwards and backwards until we reach the first entry 1. As the frieze has two rows of 
ones bounding it, we will always reach an entry 1 in each of these four directions. 
Going forwards and upwards: the first occurrence of $1$ corresponds to the diagonal $b$. Down and 
forwards: diagonal $d$. Backwards down from the entry corresponding to $1$: diagonal $c$ and backwards up: 
diagonal $e$. 
If we compare with the coordinate system for friezes of Section~\ref{ssec:fr-pattern}, the two rays 
through the object corresponding 
to diagonal $a=[kl]$ are the entries $m_{i,l}$ (with $i$ varying) and $m_{k,j}$ (with $j$ varying).

In the frieze or in the AR quiver, we give the four segments between the entry $1$ corresponding to 
$a$ and the entries corresponding to $b,c,d$ and $e$ names (see Figure~\ref{fig:regions} for a larger example containing these paths). 
Whereas $a$ is always a diagonal, $b,c,d,e$ may be boundary segments. If $b$ is a diagonal, the ray through $P_a[1]$ goes through 
$P_b[1]$, and if $b$ is a boundary segment, say $b=[i,i+1]$ (with $a=[ij]$) this ray goes through $Q_{x_i}$. By abuse of notation, it will be more convenient to 
write this projective-injective as $P_b[1]$ or as $P_{x_i}[1]$ (if we want to emphasize that it is an object 
of the Frobenius category $\mathcal C_f$ that does not live in $\mathcal C$).

Let $\mathfrak{e}$ and $\mathfrak{c}$ denote the unique sectional paths in $\mathcal{C}_f$ starting at $P_a[1]$ and ending at $P_b[1]$ and 
$P_d[1]$ respectively, but not containing $P_b[1]$ or $P_d[1]$. 
Similarly, let $\mathfrak{b}$ and $\mathfrak{d}$ denote the sectional paths in $\mathcal{C}_f$ starting at $P_e[1]$ and $P_c[1]$ respectively and 
ending at $P_a[1]$, not containing 
$P_e[1]$, $P_c[1]$, see Figure~\ref{fig:regions}.

Note that $b$ and $d$ are opposite sides of the quadrilateral determined by $a$. In particular, the corresponding diagonals do not share endpoints. 
In other words, $P_b[1]$ and $P_d[1]$ do not lie on a common ray in the AR quiver. 
So by the combinatorics of $\mathcal{C}_f$ there exist two distinct sectional paths
starting at $P_b[1], P_d[1]$. These sectional paths both go through 
$S_a$.  Let $\mathfrak{c}^a, \mathfrak{e}^a$  denote these paths starting at $P_b[1]$ and at $P_d[1]$, up to $S_a$, 
but not including $P_b [1], P_d[1]$ respectively. 
Observe that the composition of $\mathfrak{e}$ with $\mathfrak{c}^a$ and the composition of $\mathfrak{c}$ with $\mathfrak{e}^a$ are not sectional, 
see Figure~\ref{fig:regions}. 
Similarly, let $\mathfrak{d}_a,\mathfrak{b}_a$ denote the two distinct sectional paths starting at $S_a$ 
and ending at $P_e[1], P_c[1]$ respectively but not including $P_e[1], P_c[1]$.  Note that the composition of $\mathfrak{c}^a$ with 
$\mathfrak{b}_a$ and the composition of $\mathfrak{e}^a$ with $\mathfrak{d}_a$ are not sectional.

%%%%%%
\subsection{Diagonal defines subsets of indecomposables}\label{ssec:regions}
%%%%%%
$\ $

For $x$ a diagonal in the triangulation $\mathcal{T}$ and $P_x$ the corresponding projective indecomposable, 
we write $\mathcal X$ for the set of indecomposable $B$-modules having 
a non-zero homomorphism from $P_x$ into them, 
$\mathcal X=\{M \in \mbox{ind}\,B\mid \Hom(P_x,M)\ne 0\}$.  Given a $B$-module $M$, its \emph{support} is the full subquiver $\text{supp}(M)$ 
of $Q_{\mathcal{T}}$ 
generated by all vertices $x$ of $Q_{\mathcal{T}}$ such that $M\in \mathcal{X}$.  
It is well known that the support of an indecomposable module is connected.

If $x$ is a boundary segment, we set $\mathcal X$ to be the empty set (there is no projective 
indecomposable associated to $x$, so there are no indecomposables reached). 

We use the notation above to describe the regions in the frieze. 
Thus, if $x,y$ are diagonals or boundary segments, we write 
$\mathcal X\cap \mathcal Y$ for the indecomposable objects in $\mathcal{C}$ that 
have $x$ and $y$ in their support.

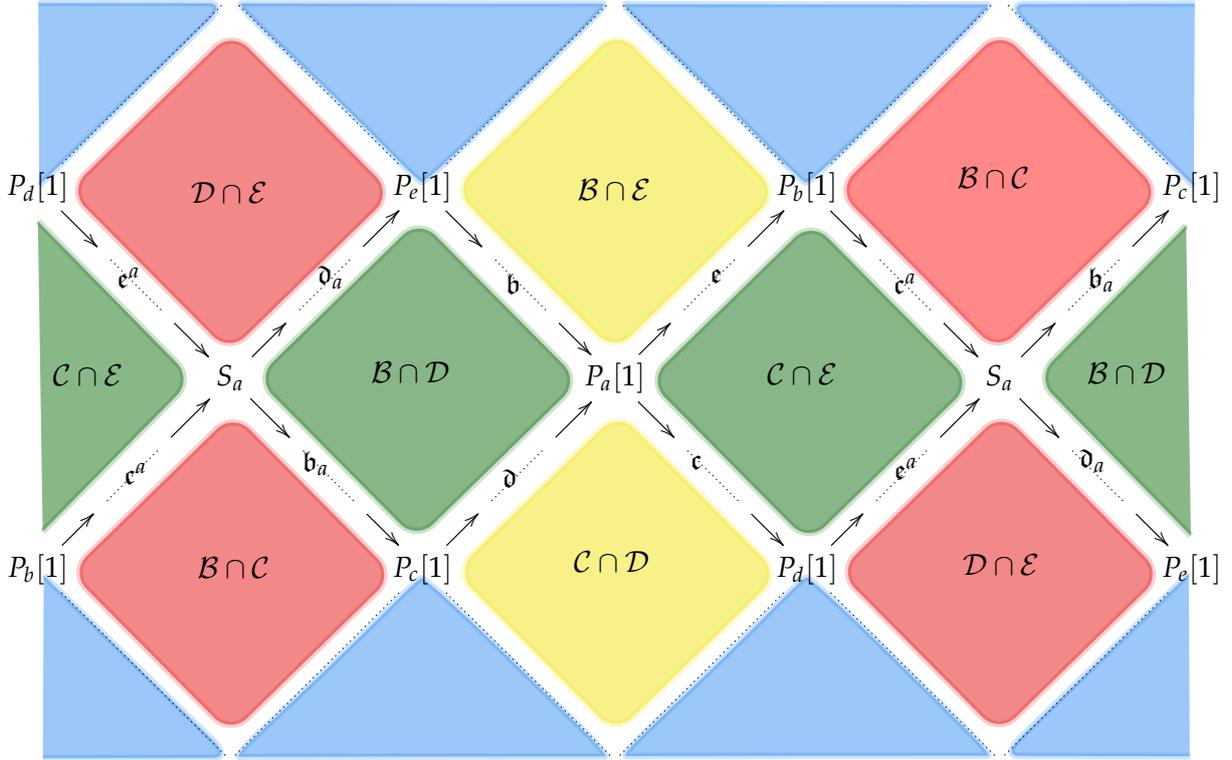
\begin{figure}
\hspace*{-2.6cm}{\begingroup%
  \makeatletter%
  \providecommand\color[2][]{%
    \errmessage{(Inkscape) Color is used for the text in Inkscape, but the package 'color.sty' is not loaded}%
    \renewcommand\color[2][]{}%
  }%
  \providecommand\transparent[1]{%
    \errmessage{(Inkscape) Transparency is used (non-zero) for the text in Inkscape, but the package 'transparent.sty' is not loaded}%
    \renewcommand\transparent[1]{}%
  }%
  \providecommand\rotatebox[2]{#2}%
  \ifx\svgwidth\undefined%
    \setlength{\unitlength}{564.02088131bp}%
    \ifx\svgscale\undefined%
      \relax%
    \else%
      \setlength{\unitlength}{\unitlength * \real{\svgscale}}%
    \fi%
  \else%
    \setlength{\unitlength}{\svgwidth}%
  \fi%
  \global\let\svgwidth\undefined%
  \global\let\svgscale\undefined%
  \makeatother%
  \begin{picture}(1,0.50983981)%
    \put(0,0){\includegraphics[width=\unitlength]{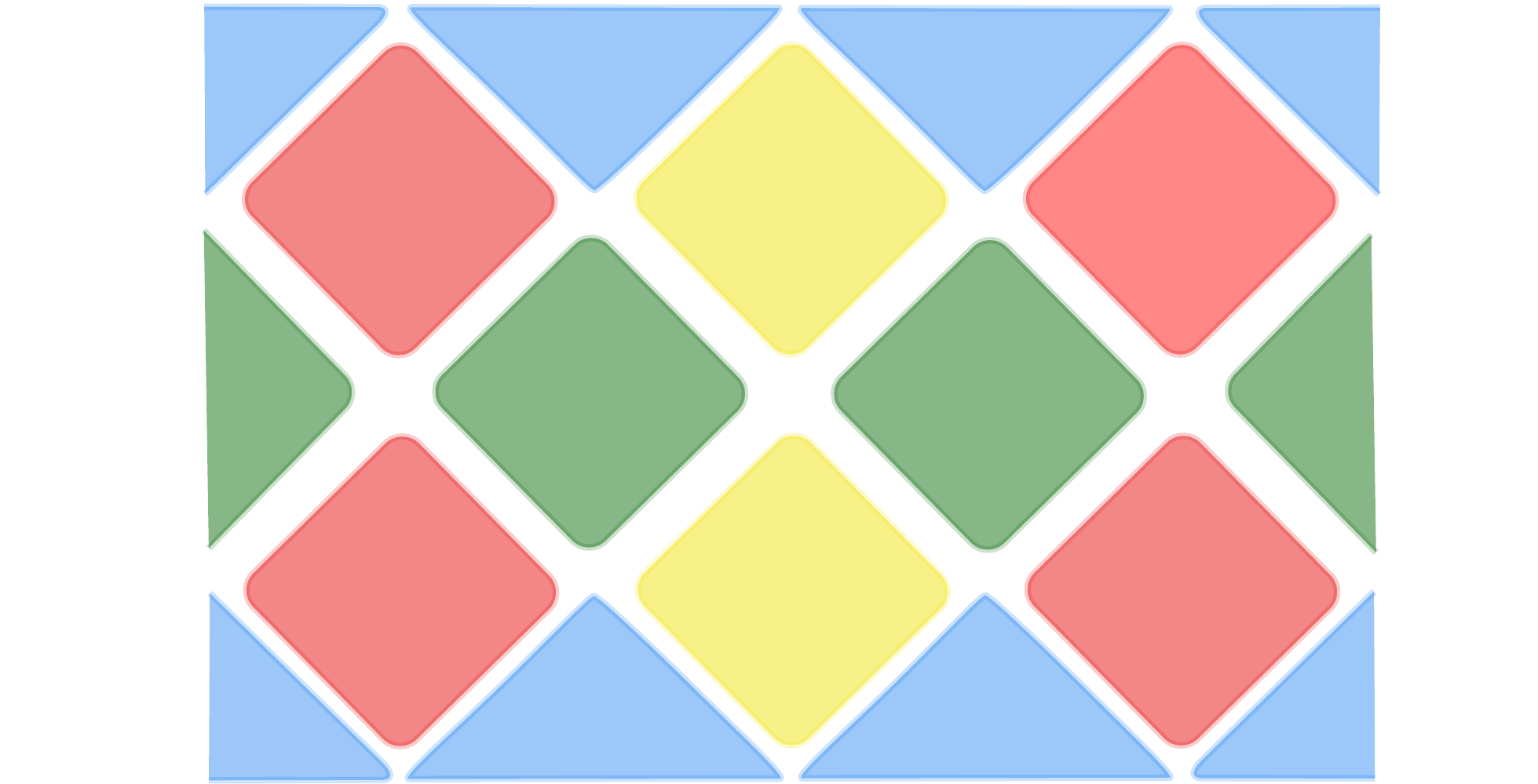}}%
    \put(-0.00030384,0.51277949){\color[rgb]{0,0,0}\makebox(0,0)[lt]{\begin{minipage}{0.61243133\unitlength}\raggedright $\xymatrix@!C=0pt@!R=0pt{&& && && \ar@{..}[dddrrr]&& && && \ar@{..}[dddrrr]&& && && \ar@{..}[dddrrr]&&&&&&\\&&&&&&&&&&&&&&&&&&&&&&&&&&&&&\\&&  &&    &&    &&  &&  &&  &&  &&  &&  &&  &&&\\&&&P_d[1] \ar@{..}[uuurrr]\ar[dr] &&  &&  && P_e[1]\ar@{..}[uuurrr] \ar[dr] && && && P_b[1]\ar@{..}[uuurrr]\ar[dr] && && && P_c[1]\\&&&&\ar@{..}[dr]&&&&\ar[ur]&&  \ar@{..}[dr]&&&&\ar[ur] &&\ar@{..}[dr]&&&&\ar[ur]&&\\&&&&&\ar[dr]&&\ar@{..}[ur]&&  &&  \ar[dr]&&  \ar@{..}[ur]&&&& \ar[dr] &&\ar@{..}[ur] &&&&\\&& && && S_a\ar[ur]\ar[dr]&&  &&  &&  P_a[1]\ar[ur]\ar[dr] &&&& && S_a \ar[ur]\ar[dr]\\&&&&&\ar[ur]&&\ar@{..}[dr]&&  &&  \ar[ur]&&  \ar@{..}[dr] && && \ar[ur] &&\ar@{..}[dr] \\&&&&\ar@{..}[ur]&&  &&  \ar[dr]&& \ar@{..}[ur]&&&& \ar[dr]&& \ar@{..}[ur]&& && \ar[dr]\\&&&P_b[1]\ar@{..}[dddrrr]\ar[ur]&&  &&  && P_c[1]\ar@{..}[dddrrr] \ar[ur]&& && && P_d[1] \ar@{..}[dddrrr]\ar[ur]&& && &&P_e[1]\\&& &&  &&  &&  &&  &&   &&   &&  &&  &&  &&\\&&&  &&  && &&  &&   && &&  &&  &&  &&\\&& && && \ar@{..}[uuurrr]&&  && && \ar@{..}[uuurrr]  &&  &&  &&  \ar@{..}[uuurrr]&&&&&&}$\\ \end{minipage}}}%
    \put(-0.09249773,0.16633676){\color[rgb]{0,0,0}\makebox(0,0)[lb]{\smash{
}}}%
    \put(0.35275657,0.26656289){\color[rgb]{0,0,0}\makebox(0,0)[lt]{\begin{minipage}{0.12737481\unitlength}\raggedright $\mathcal{B}\cap\mathcal{D}$\end{minipage}}}%
    \put(0.61469497,0.26530528){\color[rgb]{0,0,0}\makebox(0,0)[lt]{\begin{minipage}{0.12737481\unitlength}\raggedright $\mathcal{C}\cap\mathcal{E}$\end{minipage}}}%
    \put(0.48681867,0.13993638){\color[rgb]{0,0,0}\makebox(0,0)[lt]{\begin{minipage}{0.12737481\unitlength}\raggedright $\mathcal{C}\cap\mathcal{D}$\end{minipage}}}%
    \put(0.48932607,0.38766533){\color[rgb]{0,0,0}\makebox(0,0)[lt]{\begin{minipage}{0.12737481\unitlength}\raggedright $\mathcal{B}\cap\mathcal{E}$\end{minipage}}}%
    \put(0.74056533,0.39568894){\color[rgb]{0,0,0}\makebox(0,0)[lt]{\begin{minipage}{0.12737481\unitlength}\raggedright $\mathcal{B}\cap\mathcal{C}$\end{minipage}}}%
    \put(0.74407566,0.13793048){\color[rgb]{0,0,0}\makebox(0,0)[lt]{\begin{minipage}{0.12737481\unitlength}\raggedright $\mathcal{D}\cap\mathcal{E}$\end{minipage}}}%
    \put(0.23808677,0.13592457){\color[rgb]{0,0,0}\makebox(0,0)[lt]{\begin{minipage}{0.12737481\unitlength}\raggedright $\mathcal{B}\cap\mathcal{C}$\end{minipage}}}%
    \put(0.23357349,0.38565942){\color[rgb]{0,0,0}\makebox(0,0)[lt]{\begin{minipage}{0.12737481\unitlength}\raggedright $\mathcal{D}\cap\mathcal{E}$\end{minipage}}}%
    \put(0.82631766,0.26580676){\color[rgb]{0,0,0}\makebox(0,0)[lt]{\begin{minipage}{0.12737481\unitlength}\raggedright $\mathcal{B}\cap\mathcal{D}$\end{minipage}}}%
    \put(0.14230494,0.26430233){\color[rgb]{0,0,0}\makebox(0,0)[lt]{\begin{minipage}{0.12737481\unitlength}\raggedright $\mathcal{C}\cap\mathcal{E}$\end{minipage}}}%
    \put(0.44201923,0.32523551){\color[rgb]{0,0,0}\makebox(0,0)[lt]{\begin{minipage}{0.17752235\unitlength}\raggedright $\mathfrak{b}$\end{minipage}}}%
    \put(0.57732521,0.3283318){\color[rgb]{0,0,0}\makebox(0,0)[lt]{\begin{minipage}{0.17752235\unitlength}\raggedright $\mathfrak{e}$\end{minipage}}}%
    \put(0.43992088,0.1946816){\color[rgb]{0,0,0}\makebox(0,0)[lt]{\begin{minipage}{0.17752235\unitlength}\raggedright $\mathfrak{d}$\end{minipage}}}%
    \put(0.56428684,0.20370817){\color[rgb]{0,0,0}\makebox(0,0)[lt]{\begin{minipage}{0.17752235\unitlength}\raggedright $\mathfrak{c}$\end{minipage}}}%
    \put(0.82706005,0.32907707){\color[rgb]{0,0,0}\makebox(0,0)[lt]{\begin{minipage}{0.17752235\unitlength}\raggedright $\mathfrak{b}_a$\end{minipage}}}%
    \put(0.30753133,0.20671703){\color[rgb]{0,0,0}\makebox(0,0)[lt]{\begin{minipage}{0.17752235\unitlength}\raggedright $\mathfrak{b}_a$\end{minipage}}}%
    \put(0.31756084,0.33008002){\color[rgb]{0,0,0}\makebox(0,0)[lt]{\begin{minipage}{0.17752235\unitlength}\raggedright $\mathfrak{d}_a$\end{minipage}}}%
    \put(0.82140301,0.20793546){\color[rgb]{0,0,0}\makebox(0,0)[lt]{\begin{minipage}{0.17752235\unitlength}\raggedright $\mathfrak{d}_a$\end{minipage}}}%
    \put(0.18954377,0.19890891){\color[rgb]{0,0,0}\makebox(0,0)[lt]{\begin{minipage}{0.17752235\unitlength}\raggedright $\mathfrak{c}^a$\end{minipage}}}%
    \put(0.69852792,0.32528076){\color[rgb]{0,0,0}\makebox(0,0)[lt]{\begin{minipage}{0.17752235\unitlength}\raggedright $\mathfrak{c}^a$\end{minipage}}}%
    \put(0.1850169,0.32929256){\color[rgb]{0,0,0}\makebox(0,0)[lt]{\begin{minipage}{0.17752235\unitlength}\raggedright $\mathfrak{e}^a$\end{minipage}}}%
    \put(0.69889094,0.20191776){\color[rgb]{0,0,0}\makebox(0,0)[lt]{\begin{minipage}{0.17752235\unitlength}\raggedright $\mathfrak{e}^a$\end{minipage}}}%
  \end{picture}%
\endgroup%
}

  \caption{Regions in the AR quiver determined by $P_a[1]$.}
   \label{fig:regions}
\end{figure}

\begin{remark}
Let $M$ be an indecomposable $B$-module in $\mathcal{X}\cap \mathcal{Y}$ such that there exists a (unique) arrow $\alpha:\,x \to y$ in the quiver.  It follows that the right action of the element $\alpha\in B$ on $M$ is nonzero, that is $M\alpha \not=0$. 
\end{remark}

By the remark above we have the following equalities.  Note that none of the modules below are supported at $a$,  
because the same remark would imply that such modules are supported on the entire 3-cycle in $Q_{\mathcal T}$ containing $a$.  
However, this is impossible as the composition of any two arrows in a 3-cycles is zero in $B$.  
We have 

$$\mathcal{B}\cap \mathcal{E} = \{ M \in \text{ind}\,B \mid M \text{ is supported on } e\to b \}$$

$$\mathcal{C}\cap \mathcal{D} = \{ M \in \text{ind}\,B \mid M \text{ is supported on } c\to d \}$$

Moreover, since the support of an indecomposable $B$-module forms a connected subquiver of $Q$, we also have the following equalities.  

$$\mathcal{B}\cap \mathcal{C} = \{ M \in \text{ind}\,B \mid M \text{ is supported on } b\to a \to c \}$$

$$\mathcal{D}\cap \mathcal{E} = \{ M \in \text{ind}\,B \mid M \text{ is supported on } d\to a \to e \}$$

$$\mathcal{B}\cap \mathcal{D} = \{ M \in \text{ind}\,B \mid M \text{ is supported on } b\to a  \leftarrow d \}$$

$$\mathcal{C}\cap \mathcal{E} = \{ M \in \text{ind}\,B \mid M \text{ is supported on } c \leftarrow a  \to e \}$$

Finally, using similar reasoning it is easy to see that the sets described above are disjoint.  Next we describe modules lying on sectional paths defined in section \ref{ssec:diagonal-rays}. 
First, consider sectional paths starting or ending in $P_a[1]$, then we claim that  
 \[
 \mathfrak{i} = \{ M \in \text{ind}\,B \mid i \in \text{supp} (M) \subset Q_i \}\cup \{P_a[1]\}
 \]
for all $i \in \{ b, c, d, e\}$, for $Q_i$ the subquiver of $Q$ containing $i$, as in Figure~\ref{fig:quiver-divided}. 
We show that the claim holds for $i=b$, but similar arguments can be used to justify the remaining cases.  Note, that it suffices to show that a module $M \in \mathfrak{b}$ is supported on $b$ but it is not supported on $e$ or $a$.   By construction the sectional path $\mathfrak{b}$ starts at $P_e[1]$, so $0=\text{Hom}  (\tau^{-1} P_e[1], M)=\text{Hom} (P_e, M)$.  On the other hand, $\mathfrak{b}$ ends at $P_a[1]$, so $0=\text{Hom} (M, \tau P_a[1])=\text{Hom}(M, I_a)$, where $I_a$ is the injective $B$-module at $a$.  This shows that $M$ is not supported at $e$ or $a$.   
Finally, we can see from Figure~\ref{fig:regions} that $M$ has a nonzero morphism into 
$\tau P_b[1] = I_b$, provided that $b$ is not a boundary segment.  However, if $b$ is a 
boundary segment, then 
$\mathfrak{b} \cap \mathrm{Ob}(\text{mod}\,B) = \emptyset$ 
and we have $\mathfrak{b} = \{P_a[1]\}$. 
Conversely, it also follows from Figure~\ref{fig:regions} that every module $M$ supported on $b$ and some other vertices of $Q_b$ lies on $\mathfrak{b}$.  This shows the claim.   

Now consider sectional paths starting or ending in $S_a$.  Using similar arguments as above we see that

$$\mathfrak{i}^a = \{ M \in \text{ind}\,B \mid a    \in \text{supp} (M) \subset Q_i^a \}$$
for  $i \in \{c, e\}$ and 

$$\mathfrak{i}_a = \{ M \in \text{ind}\,B \mid a \in \text{supp} (M) \subset Q_i^a \} $$ 
for $i \in \{b, d\}$,  where $Q^a_i$ is the full subquiver of $Q$ on vertices of $Q_i$ and the vertex $a$.  

Finally, we define $\mathcal{F}$ to be the set of indecomposable objects of $\mathcal{C}_f$ that do not belong to 
\[
\mathcal{A}\cup \mathcal{B}\cup \mathcal{C}\cup \mathcal{D}\cup \mathcal{E} \cup \{P_a[1]\}.
\]
The region $\mathcal{F}$ is a succession of wings 
in the AR quiver of $\mathcal C_f$, with 
peaks at the $P_x[1]$ for $x\in \{b,c,d,e\}$.  That is, in the AR quiver of $\mathcal C_f$ consider two neighboured copies 
of $P_a[1]$ with the four vertices $P_b[1]$, $P_c[1]$, $P_d[1]$, $P_e[1]$. Then the indecomposables of $\mathcal{F}$ 
are the vertices in the triangular regions below these four vertices, including them (as their peaks). 
By the glide symmetry, we also have these regions at the top of the frieze. In Figure~\ref{fig:regions}, 
the wings 
are the shaded unlabelled regions at the boundary. 
It corresponds to the diagonals inside and bounding the shaded regions in Figure~\ref{fig:triangulations-1-2}.
We will see in the next section that objects in $\mathcal{F}$ do not change under mutation 
of ${T}_f$ at $P_a[1]$.  

\begin{figure}
\scalebox{.9}{\hspace{0cm}{\begingroup%
  \makeatletter%
  \providecommand\color[2][]{%
    \errmessage{(Inkscape) Color is used for the text in Inkscape, but the package 'color.sty' is not loaded}%
    \renewcommand\color[2][]{}%
  }%
  \providecommand\transparent[1]{%
    \errmessage{(Inkscape) Transparency is used (non-zero) for the text in Inkscape, but the package 'transparent.sty' is not loaded}%
    \renewcommand\transparent[1]{}%
  }%
  \providecommand\rotatebox[2]{#2}%
  \ifx\svgwidth\undefined%
    \setlength{\unitlength}{495.61414691bp}%
    \ifx\svgscale\undefined%
      \relax%
    \else%
      \setlength{\unitlength}{\unitlength * \real{\svgscale}}%
    \fi%
  \else%
    \setlength{\unitlength}{\svgwidth}%
  \fi%
  \global\let\svgwidth\undefined%
  \global\let\svgscale\undefined%
  \makeatother%
  \begin{picture}(1,0.75391361)%
    \put(0,0){\includegraphics[width=\unitlength]{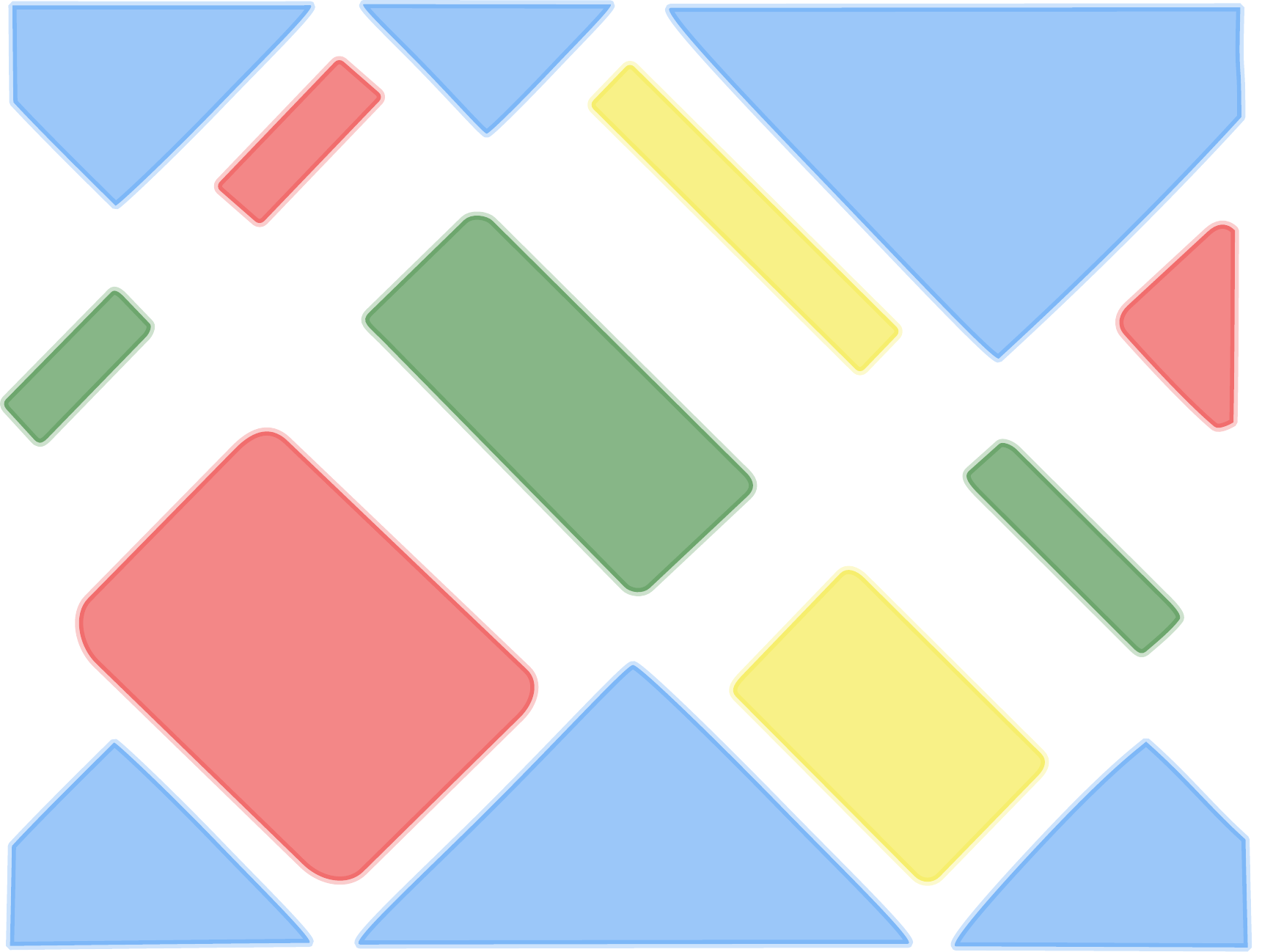}}%
    \put(0.01996563,0.74100848){\color[rgb]{0,0,0}\makebox(0,0)[lt]{\begin{minipage}{0.89375407\unitlength}\raggedright ${\xymatrix@!C=5pt@!R=5pt{& {\begin{smallmatrix}1\end{smallmatrix}} && {\begin{smallmatrix}1\end{smallmatrix}} && {\begin{smallmatrix}1\end{smallmatrix}} && {\begin{smallmatrix}1\end{smallmatrix}} && {\begin{smallmatrix}1\end{smallmatrix}} && {\begin{smallmatrix}1\end{smallmatrix}} && {\begin{smallmatrix}1\end{smallmatrix}} && {\begin{smallmatrix}1\end{smallmatrix}}  \\{\begin{smallmatrix}1\end{smallmatrix}} && {\begin{smallmatrix}2\end{smallmatrix}} && {\begin{smallmatrix}3 \\ \textcolor{red}{4}\end{smallmatrix}} && {\begin{smallmatrix}1\end{smallmatrix}} && {\begin{smallmatrix}5\\ \textcolor{red}{4}\end{smallmatrix}}  &&{\begin{smallmatrix}1 \end{smallmatrix}} && {\begin{smallmatrix}4\end{smallmatrix}} && {\begin{smallmatrix}1\end{smallmatrix}} && {\begin{smallmatrix}2\end{smallmatrix}} \\ & {\begin{smallmatrix}1 \end{smallmatrix}} && {\begin{smallmatrix}5\\ \textcolor{red}{7}\end{smallmatrix}} && {\begin{smallmatrix} 2 \\ \textcolor{red}{3} \end{smallmatrix}} && {\begin{smallmatrix}4\\ \textcolor{red}{3}\end{smallmatrix}} && {\begin{smallmatrix}  4 \\\textcolor{red}{3} \end{smallmatrix}} && {\begin{smallmatrix} 3 \end{smallmatrix}} && {\begin{smallmatrix}3 \end{smallmatrix}} && {\begin{smallmatrix}1\end{smallmatrix}} \\ {\begin{smallmatrix}4 \\ \textcolor{red}{3}\end{smallmatrix}} && {\begin{smallmatrix}2 \\ \textcolor{red}{3} \end{smallmatrix}} && {\begin{smallmatrix}3 \\ \textcolor{red}{5}\end{smallmatrix}} && {\begin{smallmatrix}7 \\ \textcolor{red}{8}\end{smallmatrix}} && {\begin{smallmatrix}3\\\textcolor{red}{2}\end{smallmatrix}} && {\begin{smallmatrix}11 \\ \textcolor{red}{8}\end{smallmatrix}} && {\begin{smallmatrix}2 \end{smallmatrix}} && {\begin{smallmatrix}2\end{smallmatrix}} && {\begin{smallmatrix}4\\ \textcolor{red}{5} \end{smallmatrix}} \\& {\begin{smallmatrix}7 \\ \textcolor{red}{8}\end{smallmatrix}} &&  {\begin{smallmatrix}1 \\ \textcolor{red}{2}\end{smallmatrix}} && {\begin{smallmatrix}10\\ \textcolor{red}{13}\end{smallmatrix}} && {\begin{smallmatrix}5 \\ \textcolor{red}{5}\end{smallmatrix}} && {\begin{smallmatrix}8 \\ \textcolor{red}{5}\end{smallmatrix}} && {\begin{smallmatrix}7\\ \textcolor{red}{5}\end{smallmatrix}} && {\begin{smallmatrix}1 \end{smallmatrix}} &&  {\begin{smallmatrix}7 \\ \textcolor{red}{9}\end{smallmatrix}}  \\{\begin{smallmatrix}5 \\ \textcolor{red}{5}\end{smallmatrix}} &&  {\begin{smallmatrix}3 \\ \textcolor{red}{5}\end{smallmatrix}} && {\begin{smallmatrix}3 \\ \textcolor{red}{5}\end{smallmatrix}} && {\begin{smallmatrix}7 \\ \textcolor{red}{8}\end{smallmatrix}} && {\begin{smallmatrix}13 \\ \textcolor{red}{12}\end{smallmatrix}} && {\begin{smallmatrix}5 \\ \textcolor{red}{3}\end{smallmatrix}} && {\begin{smallmatrix}3 \\ \textcolor{red}{2}\end{smallmatrix}} &&  {\begin{smallmatrix}3 \\ \textcolor{red}{4}\end{smallmatrix}} && {\begin{smallmatrix}5 \\ \textcolor{red}{7}\end{smallmatrix}}  \\&  {\begin{smallmatrix}2 \\ \textcolor{red}{3}\end{smallmatrix}} && {\begin{smallmatrix}8 \\ \textcolor{red}{12}\end{smallmatrix}} && {\begin{smallmatrix}2 \\ \textcolor{red}{3}\end{smallmatrix}} && {\begin{smallmatrix}18 \\ \textcolor{red}{19}\end{smallmatrix}} && {\begin{smallmatrix}8 \\ \textcolor{red}{7}\end{smallmatrix}} && {\begin{smallmatrix}2 \\ \textcolor{red}{1}\end{smallmatrix}} && {\begin{smallmatrix}8 \\ \textcolor{red}{7}\end{smallmatrix}} && {\begin{smallmatrix}2 \\ \textcolor{red}{3}\end{smallmatrix}}  \\{\begin{smallmatrix}3 \\ \textcolor{red}{4}\end{smallmatrix}} && {\begin{smallmatrix}5\\ \textcolor{red}{7}\end{smallmatrix}} && {\begin{smallmatrix}5 \\ \textcolor{red}{7}\end{smallmatrix}} && {\begin{smallmatrix}5 \\ \textcolor{red}{7}\end{smallmatrix}} && {\begin{smallmatrix}11 \\ \textcolor{red}{11}\end{smallmatrix}} && {\begin{smallmatrix}3 \\ \textcolor{red}{2}\end{smallmatrix}} && {\begin{smallmatrix}5 \\ \textcolor{red}{3}\end{smallmatrix}} &&  {\begin{smallmatrix}5 \\ \textcolor{red}{5}\end{smallmatrix}} && {\begin{smallmatrix}3 \\ \textcolor{red}{5}\end{smallmatrix}}  \\&  {\begin{smallmatrix}7\\ \textcolor{red}{9}\end{smallmatrix}} &&  {\begin{smallmatrix}3 \\ \textcolor{red}{4}\end{smallmatrix}} && {\begin{smallmatrix}12 \\ \textcolor{red}{16} \end{smallmatrix}} && {\begin{smallmatrix}3 \\ \textcolor{red}{4}\end{smallmatrix}} && {\begin{smallmatrix}4 \\ \textcolor{red}{3}\end{smallmatrix}} && {\begin{smallmatrix}7 \\ \textcolor{red}{5}\end{smallmatrix}} && {\begin{smallmatrix}3\\ \textcolor{red}{2}\end{smallmatrix}} &&  {\begin{smallmatrix}7 \\ \textcolor{red}{8}\end{smallmatrix}}  \\{\begin{smallmatrix}2\end{smallmatrix}} &&  {\begin{smallmatrix}  4\\ \textcolor{red}{5}\end{smallmatrix}} && {\begin{smallmatrix}7 \\ \textcolor{red}{9}\end{smallmatrix}} &&  {\begin{smallmatrix}7\\ \textcolor{red}{9} \end{smallmatrix}} && {\begin{smallmatrix}1 \end{smallmatrix}} && {\begin{smallmatrix}9 \\ \textcolor{red}{7}\end{smallmatrix}} && {\begin{smallmatrix}4 \\ \textcolor{red}{3}\end{smallmatrix}} &&  {\begin{smallmatrix}4\\ \textcolor{red}{3}\end{smallmatrix}} && {\begin{smallmatrix}2\\ \textcolor{red}{3}\end{smallmatrix}}  \\&  {\begin{smallmatrix}1\end{smallmatrix}} &&  {\begin{smallmatrix} 9\\ \textcolor{red}{11} \end{smallmatrix}} && {\begin{smallmatrix} 4\\\textcolor{red}{5} \end{smallmatrix}} &&  {\begin{smallmatrix}2\end{smallmatrix}} && {\begin{smallmatrix}2 \end{smallmatrix}} &&  {\begin{smallmatrix}5\\ \textcolor{red}{4}\end{smallmatrix}} &&     {\begin{smallmatrix}5\\ \textcolor{red}{4}\end{smallmatrix}} &&  {\begin{smallmatrix}1\end{smallmatrix}}  \\{\begin{smallmatrix}1\end{smallmatrix}} &&  {\begin{smallmatrix}2 \end{smallmatrix}} &&     {\begin{smallmatrix}5\\ \textcolor{red}{6}\end{smallmatrix}} &&  {\begin{smallmatrix}1\end{smallmatrix}} &&     {\begin{smallmatrix}3 \end{smallmatrix}} &&  {\begin{smallmatrix}1\end{smallmatrix}} &&     {\begin{smallmatrix}6 \\ \textcolor{red}{5}\end{smallmatrix}} &&  {\begin{smallmatrix}1\end{smallmatrix}} &&    {\begin{smallmatrix}2\end{smallmatrix}} &&\\&{\begin{smallmatrix}1\end{smallmatrix}} &&  {\begin{smallmatrix}1\end{smallmatrix}} &&     {\begin{smallmatrix}1\end{smallmatrix}} &&  {\begin{smallmatrix}1\end{smallmatrix}} &&     {\begin{smallmatrix}1\end{smallmatrix}} &&  {\begin{smallmatrix}1\end{smallmatrix}} &&     {\begin{smallmatrix}1\end{smallmatrix}} &&  {\begin{smallmatrix}1\end{smallmatrix}} &&    {\begin{smallmatrix}1\end{smallmatrix}} &&}}$\end{minipage}}}%
  \end{picture}%
\endgroup%
}}

  \caption{Frieze pattern of Example~\ref{ex:quiver-triangul-flip}. Red entries: after flip of diagonal 1}
   \label{fig:frieze-entries}
\end{figure}

\begin{ex}\label{ex:quiver-triangul-flip}
We consider the triangulation $\mathcal{T}$ of a 14-gon, see left hand of Figure~\ref{fig:triangulations-1-2} 
and the triangulation $\mathcal{T}'=\mu_1(\mathcal{T})$ obtained by flipping diagonal $1$.  

The quivers of $\mathcal{T}$ and of $\mathcal{T}'$ are given below. Note that the quiver $Q$ is the same as in Example~\ref{E:category_frieze}.

$$\xymatrix@C=10pt@R=10pt{&&8\ar[dr]&&&&9\ar[dr] &&    &&& &&8\ar[dr]&&&&9\ar[dr] &&\\
Q: &7\ar[ur]&&3\ar[ll]\ar[dr]&&2\ar[ll]\ar[ur]&&10\ar[ll]\ar[d]&& &&    Q': &7\ar[ur]&&3\ar[dd]\ar[ll]&&2 \ar[ur]\ar[dl]&&10\ar[ll]\ar[d]\\
&&&&1\ar[ur]\ar[dl] &&&11&& &&&&  &&1\ar[ul]\ar[dr] &&&11&&\\
&6&&5\ar[ll]\ar[rr]&&4\ar[ul] && && && &6&&5\ar[ll]\ar[ur]&&4\ar[uu]}$$

Figure~\ref{fig:ARquiver} shows the Auslander-Reiten quiver of the cluster category 
$\mathcal C_f$ for $Q$. 

In Figure~\ref{fig:frieze-entries} (Section~\ref{sec:mutating}), the frieze patterns of $T$ and of $T'$ are given. 

\begin{figure}

\scalebox{.7}{\begingroup%
  \makeatletter%
  \providecommand\color[2][]{%
    \errmessage{(Inkscape) Color is used for the text in Inkscape, but the package 'color.sty' is not loaded}%
    \renewcommand\color[2][]{}%
  }%
  \providecommand\transparent[1]{%
    \errmessage{(Inkscape) Transparency is used (non-zero) for the text in Inkscape, but the package 'transparent.sty' is not loaded}%
    \renewcommand\transparent[1]{}%
  }%
  \providecommand\rotatebox[2]{#2}%
  \ifx\svgwidth\undefined%
    \setlength{\unitlength}{616.18436966bp}%
    \ifx\svgscale\undefined%
      \relax%
    \else%
      \setlength{\unitlength}{\unitlength * \real{\svgscale}}%
    \fi%
  \else%
    \setlength{\unitlength}{\svgwidth}%
  \fi%
  \global\let\svgwidth\undefined%
  \global\let\svgscale\undefined%
  \makeatother%
  \begin{picture}(1,0.42069957)%
    \put(0,0){\includegraphics[width=\unitlength]{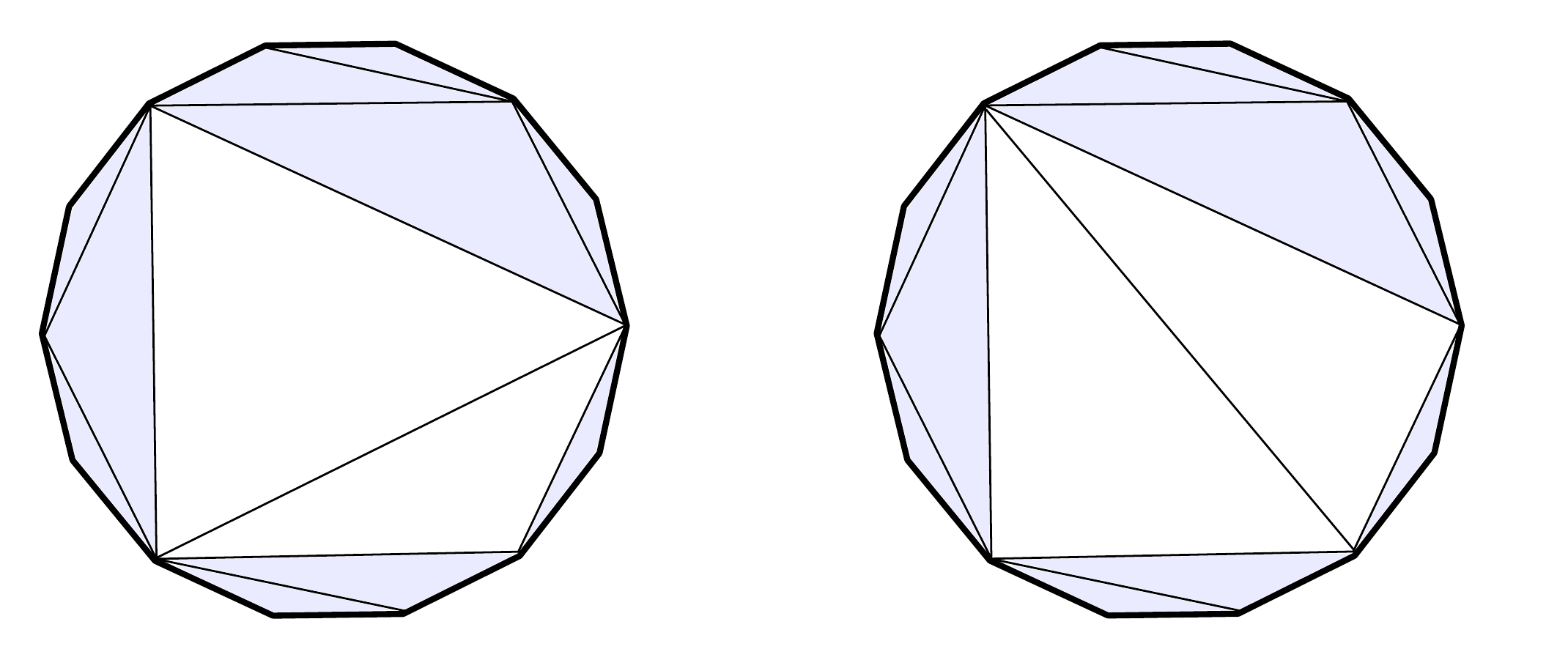}}%
    \put(0.24932661,0.17096629){\color[rgb]{0,0,0}\makebox(0,0)[lt]{\begin{minipage}{0.04196782\unitlength}\raggedright 1\end{minipage}}}%
    \put(0.24211336,0.25883638){\color[rgb]{0,0,0}\makebox(0,0)[lb]{\smash{2}}}%
    \put(0.10243928,0.17293351){\color[rgb]{0,0,0}\makebox(0,0)[lb]{\smash{3}}}%
    \put(0.34244269,0.12703123){\color[rgb]{0,0,0}\makebox(0,0)[lb]{\smash{4}}}%
    \put(0.27227776,0.07194847){\color[rgb]{0,0,0}\makebox(0,0)[lb]{\smash{5}}}%
    \put(0.23949037,0.03784965){\color[rgb]{0,0,0}\makebox(0,0)[lb]{\smash{6}}}%
    \put(0.06047147,0.25818064){\color[rgb]{0,0,0}\makebox(0,0)[lb]{\smash{7}}}%
    \put(0.05784847,0.15391686){\color[rgb]{0,0,0}\makebox(0,0)[lb]{\smash{8}}}%
    \put(0.34309838,0.27982027){\color[rgb]{0,0,0}\makebox(0,0)[lb]{\smash{9}}}%
    \put(0.22637544,0.33555877){\color[rgb]{0,0,0}\makebox(0,0)[lb]{\smash{10}}}%
    \put(0.16670246,0.36965763){\color[rgb]{0,0,0}\makebox(0,0)[lb]{\smash{11}}}%
    \put(0.04479613,0.32484481){\color[rgb]{0,0,0}\makebox(0,0)[lb]{\smash{$x_1$}}}%
    \put(0.10136548,0.38141411){\color[rgb]{0,0,0}\makebox(0,0)[lb]{\smash{$x_2$
}}}%
    \put(0.19688421,0.40459828){\color[rgb]{0,0,0}\makebox(0,0)[lb]{\smash{$x_3$}}}%
    \put(0.29240291,0.38512361){\color[rgb]{0,0,0}\makebox(0,0)[lb]{\smash{$x_4$}}}%
    \put(0.36659218,0.33040898){\color[rgb]{0,0,0}\makebox(0,0)[lb]{\smash{$x_5$}}}%
    \put(0.40368685,0.24694601){\color[rgb]{0,0,0}\makebox(0,0)[lb]{\smash{$x_6$}}}%
    \put(0.40183217,0.1486452){\color[rgb]{0,0,0}\makebox(0,0)[lb]{\smash{$x_7$}}}%
    \put(0.37308379,0.08372956){\color[rgb]{0,0,0}\makebox(0,0)[lb]{\smash{$x_8$}}}%
    \put(0.30353132,0.02901496){\color[rgb]{0,0,0}\makebox(0,0)[lb]{\smash{$x_9$}}}%
    \put(0.19410208,0.00397608){\color[rgb]{0,0,0}\makebox(0,0)[lb]{\smash{$x_{10}$}}}%
    \put(0.10600228,0.02530549){\color[rgb]{0,0,0}\makebox(0,0)[lb]{\smash{$x_{11}$}}}%
    \put(0.03737723,0.08187485){\color[rgb]{0,0,0}\makebox(0,0)[lb]{\smash{$x_{12}$}}}%
    \put(0.00677412,0.15420942){\color[rgb]{0,0,0}\makebox(0,0)[lb]{\smash{$x_{13}$}}}%
    \put(-0.00157218,0.24880072){\color[rgb]{0,0,0}\makebox(0,0)[lb]{\smash{$x_{14}$}}}%
    \put(0.77077233,0.16737265){\color[rgb]{0,0,0}\makebox(0,0)[lt]{\begin{minipage}{0.04196782\unitlength}\raggedright 1\end{minipage}}}%
    \put(0.77468751,0.25895219){\color[rgb]{0,0,0}\makebox(0,0)[lb]{\smash{2}}}%
    \put(0.63501338,0.17304932){\color[rgb]{0,0,0}\makebox(0,0)[lb]{\smash{3}}}%
    \put(0.8750168,0.12714704){\color[rgb]{0,0,0}\makebox(0,0)[lb]{\smash{4}}}%
    \put(0.80485189,0.07206428){\color[rgb]{0,0,0}\makebox(0,0)[lb]{\smash{5}}}%
    \put(0.7720645,0.03796546){\color[rgb]{0,0,0}\makebox(0,0)[lb]{\smash{6}}}%
    \put(0.59304557,0.25829646){\color[rgb]{0,0,0}\makebox(0,0)[lb]{\smash{7}}}%
    \put(0.59042259,0.15403267){\color[rgb]{0,0,0}\makebox(0,0)[lb]{\smash{8}}}%
    \put(0.87567253,0.27993608){\color[rgb]{0,0,0}\makebox(0,0)[lb]{\smash{9}}}%
    \put(0.75894959,0.33567459){\color[rgb]{0,0,0}\makebox(0,0)[lb]{\smash{10}}}%
    \put(0.69927657,0.36977344){\color[rgb]{0,0,0}\makebox(0,0)[lb]{\smash{11}}}%
    \put(0.57737025,0.32496062){\color[rgb]{0,0,0}\makebox(0,0)[lb]{\smash{$x_1$}}}%
    \put(0.63393961,0.38152993){\color[rgb]{0,0,0}\makebox(0,0)[lb]{\smash{$x_2$
}}}%
    \put(0.7294583,0.4047141){\color[rgb]{0,0,0}\makebox(0,0)[lb]{\smash{$x_3$}}}%
    \put(0.824977,0.38523943){\color[rgb]{0,0,0}\makebox(0,0)[lb]{\smash{$x_4$}}}%
    \put(0.89916633,0.33052479){\color[rgb]{0,0,0}\makebox(0,0)[lb]{\smash{$x_5$}}}%
    \put(0.93626096,0.24706182){\color[rgb]{0,0,0}\makebox(0,0)[lb]{\smash{$x_6$}}}%
    \put(0.93440628,0.14876101){\color[rgb]{0,0,0}\makebox(0,0)[lb]{\smash{$x_7$}}}%
    \put(0.9056579,0.08384537){\color[rgb]{0,0,0}\makebox(0,0)[lb]{\smash{$x_8$}}}%
    \put(0.83610545,0.02913077){\color[rgb]{0,0,0}\makebox(0,0)[lb]{\smash{$x_9$}}}%
    \put(0.72667617,0.0040919){\color[rgb]{0,0,0}\makebox(0,0)[lb]{\smash{$x_{10}$}}}%
    \put(0.63857641,0.0254213){\color[rgb]{0,0,0}\makebox(0,0)[lb]{\smash{$x_{11}$}}}%
    \put(0.56995135,0.08199066){\color[rgb]{0,0,0}\makebox(0,0)[lb]{\smash{$x_{12}$}}}%
    \put(0.53934825,0.15432524){\color[rgb]{0,0,0}\makebox(0,0)[lb]{\smash{$x_{13}$}}}%
    \put(0.53100193,0.24891654){\color[rgb]{0,0,0}\makebox(0,0)[lb]{\smash{$x_{14}$}}}%
  \end{picture}%
\endgroup%
}

\caption{Triangulations $\mathcal{T}$ and $\mathcal{T}'=\mu_1(\mathcal{T})$}
\label{fig:triangulations-1-2}
\end{figure}
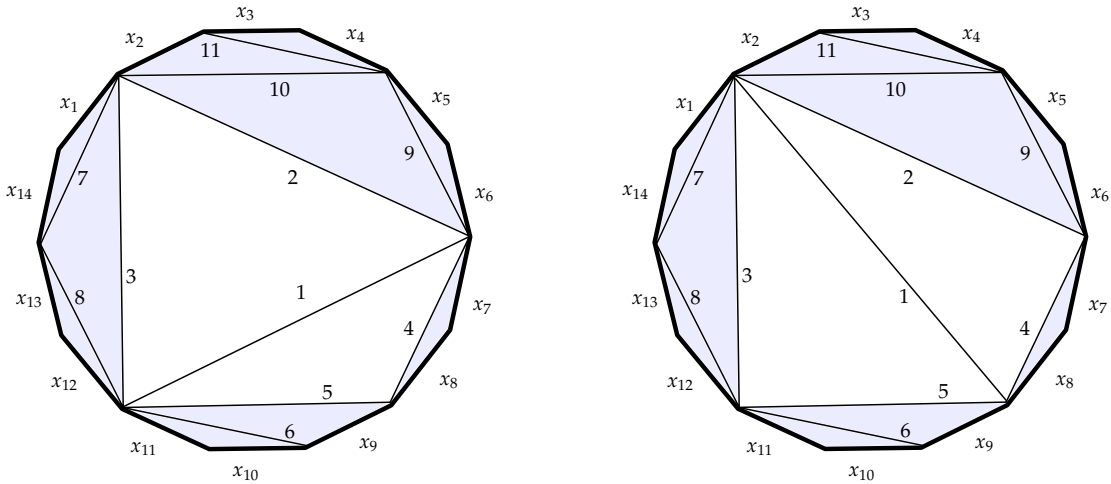

\begin{figure}

\scalebox{0.9}{\hspace*{0cm}{\begingroup%
  \makeatletter%
  \providecommand\color[2][]{%
    \errmessage{(Inkscape) Color is used for the text in Inkscape, but the package 'color.sty' is not loaded}%
    \renewcommand\color[2][]{}%
  }%
  \providecommand\transparent[1]{%
    \errmessage{(Inkscape) Transparency is used (non-zero) for the text in Inkscape, but the package 'transparent.sty' is not loaded}%
    \renewcommand\transparent[1]{}%
  }%
  \providecommand\rotatebox[2]{#2}%
  \ifx\svgwidth\undefined%
    \setlength{\unitlength}{516.69798254bp}%
    \ifx\svgscale\undefined%
      \relax%
    \else%
      \setlength{\unitlength}{\unitlength * \real{\svgscale}}%
    \fi%
  \else%
    \setlength{\unitlength}{\svgwidth}%
  \fi%
  \global\let\svgwidth\undefined%
  \global\let\svgscale\undefined%
  \makeatother%
  \begin{picture}(1,0.72073294)%
    \put(0,0){\includegraphics[width=\unitlength]{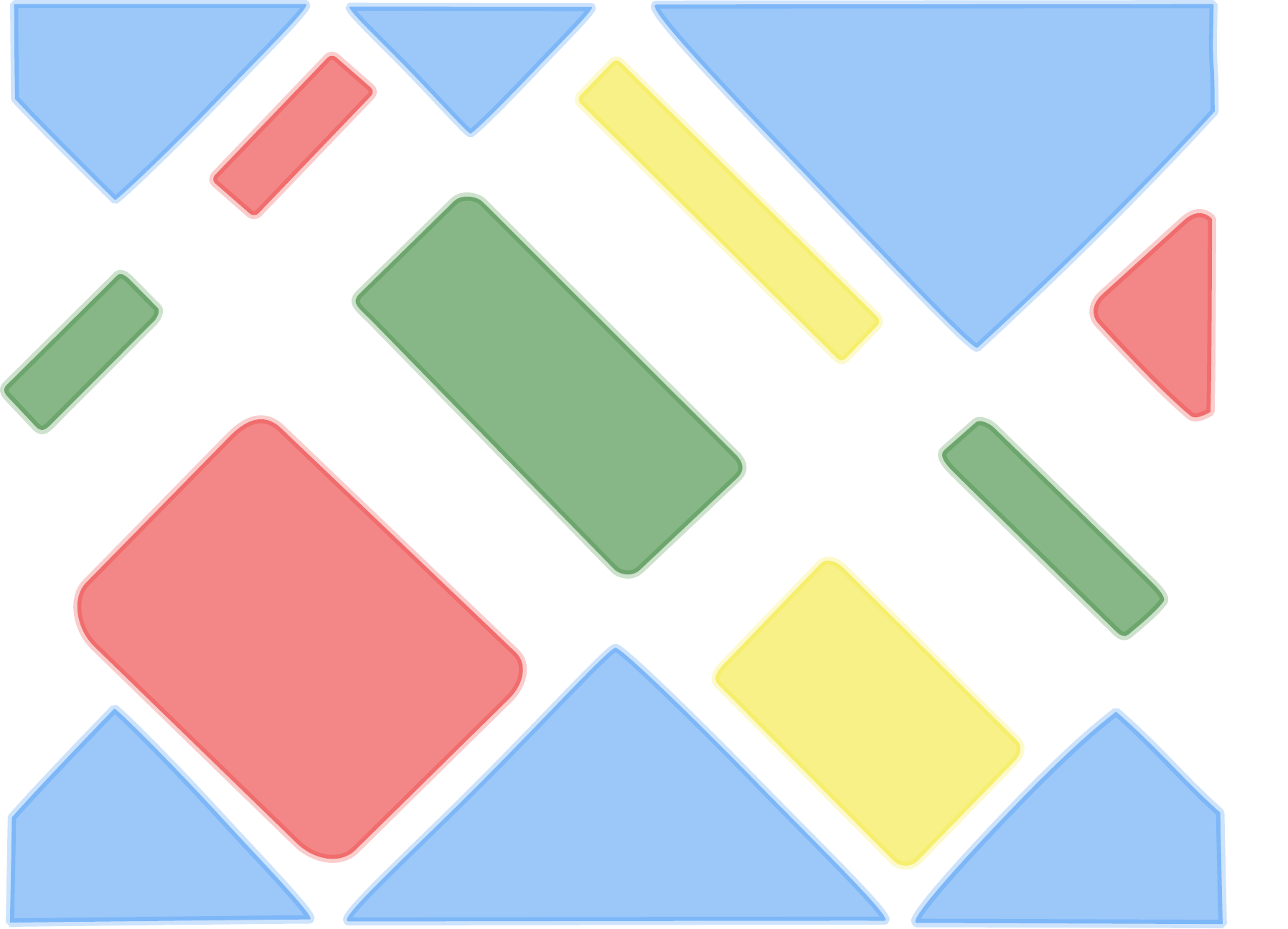}}%
    \put(0.00272146,0.71144678){\color[rgb]{0,0,0}\makebox(0,0)[lt]{\begin{minipage}{0.99987766\unitlength}\raggedright $\xymatrix@!C=5pt@!R=5pt{&{\begin{smallmatrix}P_{x_9}[1] \end{smallmatrix}}\ar[dr]&&{\begin{smallmatrix}P_{x_{8}}[1]\end{smallmatrix}}\ar[dr]&&{\begin{smallmatrix}P_{x_{7}}[1]\end{smallmatrix}}\ar[dr]&&{\begin{smallmatrix}P_{x_{6}}[1]\end{smallmatrix}}\ar[dr]&&{\begin{smallmatrix}P_{x_{5}}[1]\end{smallmatrix}}\ar[dr]&&{\begin{smallmatrix}P_{x_{4}}[1]\end{smallmatrix}}\ar[dr]&&{\begin{smallmatrix}P_{x_{3}}[1]\end{smallmatrix}}\ar[dr]&&{\begin{smallmatrix}P_{x_{2}}[1]\end{smallmatrix}}\ar[dr]\\{\begin{smallmatrix}P_6[1]\end{smallmatrix}}\ar[dr]\ar[ur]&&{\begin{smallmatrix}6\end{smallmatrix}}\ar[dr]\ar[ur]&&{\begin{smallmatrix}5\\4\end{smallmatrix}}\ar[dr]\ar[ur]&&{\begin{smallmatrix}P_4[1]\end{smallmatrix}}\ar[dr]\ar[ur]&&{\begin{smallmatrix}4\\1\\2\\9\end{smallmatrix}}\ar[dr]\ar[ur]&&{\begin{smallmatrix}P_9[1]\end{smallmatrix}}\ar[dr]\ar[ur]&&{\begin{smallmatrix}9\\10\\11\end{smallmatrix}}\ar[dr]\ar[ur]&&{\begin{smallmatrix}P_{11}[1]\end{smallmatrix}}\ar[dr]\ar[ur]&&{\begin{smallmatrix}11\end{smallmatrix}}\\&{\begin{smallmatrix}P_5[1]\end{smallmatrix}}\ar[ur]\ar[dr]&&{\begin{smallmatrix}5\\6\;4\end{smallmatrix}}\ar[ur]\ar[dr]&&{\begin{smallmatrix}5\end{smallmatrix}}\ar[ur]\ar[dr]&&{\begin{smallmatrix}1\\2\\9\end{smallmatrix}}\ar[ur]\ar[dr]&&{\begin{smallmatrix}4\\1\\2\end{smallmatrix}}\ar[ur]\ar[dr]&&{\begin{smallmatrix}10\\11\end{smallmatrix}}\ar[ur]\ar[dr]&&{\begin{smallmatrix}9\\10\end{smallmatrix}}\ar[ur]\ar[dr]&&{\begin{smallmatrix}P_{10}[1]\end{smallmatrix}}\ar[ur]\ar[dr]\\{\begin{smallmatrix}8\\3\\1\end{smallmatrix}}\ar[ur]\ar[dr]&&{\begin{smallmatrix}4\end{smallmatrix}}\ar[ur]\ar[dr]&&{\begin{smallmatrix}5\\6\end{smallmatrix}}\ar[ur]\ar[dr]&&{\begin{smallmatrix}1\\2\;5\\9\;\;\;\end{smallmatrix}}\ar[ur]\ar[dr]&&{\begin{smallmatrix}1\\2\end{smallmatrix}}\ar[ur]\ar[dr]&&{\begin{smallmatrix}\;\;\;\;\;\;\;4\\\;\;\;10\;1\\11\;2\end{smallmatrix}}\ar[ur]\ar[dr]&&{\begin{smallmatrix}10\end{smallmatrix}}\ar[ur]\ar[dr]&&{\begin{smallmatrix}9\end{smallmatrix}}\ar[ur]\ar[dr]&&{\begin{smallmatrix}2\\3\\7\end{smallmatrix}}\\&{\begin{smallmatrix}\;\;\;\;8\\4\;3\\1\end{smallmatrix}}\ar[ur]\ar[dr]&&{\begin{smallmatrix}P_1[1]\end{smallmatrix}}\ar[ur]\ar[dr]&&{\begin{smallmatrix}1\\2\;5\\9\;6\end{smallmatrix}}\ar[ur]\ar[dr]&&{\begin{smallmatrix}1\\2\;5\end{smallmatrix}}\ar[ur]\ar[dr]&&{\begin{smallmatrix}\;\;\;10\;\;1\\11\;2\end{smallmatrix}}\ar[ur]\ar[dr]&&{\begin{smallmatrix}\;\;\;\;4\\10\;1\\\;2\end{smallmatrix}}\ar[ur]\ar[dr]&&{\begin{smallmatrix}P_2[1]\end{smallmatrix}}\ar[ur]\ar[dr]&&{\begin{smallmatrix}2\\9\;3\\\;\;\;\;7\end{smallmatrix}}\ar[ur]\ar[dr]\\{\begin{smallmatrix}4\;3\\1\end{smallmatrix}}\ar[ur]\ar[dr]&&{\begin{smallmatrix}8\\3\end{smallmatrix}}\ar[ur]\ar[dr]&&{\begin{smallmatrix}2\\9\end{smallmatrix}}\ar[ur]\ar[dr]&&{\begin{smallmatrix}1\\2\;5\\\;\;\;6\end{smallmatrix}}\ar[ur]\ar[dr]&&{\begin{smallmatrix}\;\;10\;1\\11\;2\;5\end{smallmatrix}}\ar[ur]\ar[dr]&&{\begin{smallmatrix}10\;1\\2\end{smallmatrix}}\ar[ur]\ar[dr]&&{\begin{smallmatrix}4\\1\end{smallmatrix}}\ar[ur]\ar[dr]&&{\begin{smallmatrix}3\\7\end{smallmatrix}}\ar[ur]\ar[dr]&&{\begin{smallmatrix}2\\9\;3\end{smallmatrix}}\\&{\begin{smallmatrix}3\end{smallmatrix}}\ar[ur]\ar[dr]&&{\begin{smallmatrix}\;\;2\;8\\9\;3\end{smallmatrix}}\ar[ur]\ar[dr]&&{\begin{smallmatrix}2\end{smallmatrix}}\ar[ur]\ar[dr]&&{\begin{smallmatrix}\;\;10\;1\\11\;2\;5\\\;\;\;\;\;\;\;\;6\end{smallmatrix}}\ar[ur]\ar[dr]&&{\begin{smallmatrix}10\;1\\2\;5\end{smallmatrix}}\ar[ur]\ar[dr]&&{\begin{smallmatrix}1\end{smallmatrix}}\ar[ur]\ar[dr]&&{\begin{smallmatrix}4\;3\\\;\;1\;7\end{smallmatrix}}\ar[ur]\ar[dr]&&{\begin{smallmatrix}3\end{smallmatrix}}\ar[ur]\ar[dr]\\{\begin{smallmatrix}3\\7\end{smallmatrix}}\ar[ur]\ar[dr]&&{\begin{smallmatrix}2\\9\;3\end{smallmatrix}}\ar[ur]\ar[dr]&&{\begin{smallmatrix}8\;2\\3\end{smallmatrix}}\ar[ur]\ar[dr]&&{\begin{smallmatrix}10\\2\;11\end{smallmatrix}}\ar[ur]\ar[dr]&&{\begin{smallmatrix}10\;\;1\\\;\;\;2\;5\\\;\;\;\;\;\;\;6\end{smallmatrix}}\ar[ur]\ar[dr]&&{\begin{smallmatrix}1\\5\end{smallmatrix}}\ar[ur]\ar[dr]&&{\begin{smallmatrix}3\\1\;7\end{smallmatrix}}\ar[ur]\ar[dr]&&{\begin{smallmatrix}4\;3\\1\end{smallmatrix}}\ar[ur]\ar[dr]&&{\begin{smallmatrix}8\\3\end{smallmatrix}}\\&{\begin{smallmatrix}2\\9\;3\\\;\;\;\;7\end{smallmatrix}}\ar[ur]\ar[dr]&&{\begin{smallmatrix}2\\3\end{smallmatrix}}\ar[ur]\ar[dr]&&{\begin{smallmatrix}\;\;\;10\\8\;2\;11\\3\;\;\;\;\end{smallmatrix}}\ar[ur]\ar[dr]&&{\begin{smallmatrix}10\\2\end{smallmatrix}}\ar[ur]\ar[dr]&&{\begin{smallmatrix}1\\5\\6\end{smallmatrix}}\ar[ur]\ar[dr]&&{\begin{smallmatrix}3\\1\;7\\5\;\;\;\;\end{smallmatrix}}\ar[ur]\ar[dr]&&{\begin{smallmatrix}3\\1\end{smallmatrix}}\ar[ur]\ar[dr]&&{\begin{smallmatrix}\;\;\;\;8\\4\;3\\1\end{smallmatrix}}\ar[ur]\ar[dr]\\{\begin{smallmatrix}9\end{smallmatrix}}\ar[ur]\ar[dr]&&{\begin{smallmatrix}2\\3\\7\end{smallmatrix}}\ar[ur]\ar[dr]&&{\begin{smallmatrix}10\\2\;11\\3\;\;\;\;\;\;\end{smallmatrix}}\ar[ur]\ar[dr]&&{\begin{smallmatrix}\;\;\;\;10\\8\;2\\3\;\end{smallmatrix}}\ar[ur]\ar[dr]&&{\begin{smallmatrix}P_3[1]\end{smallmatrix}}\ar[ur]\ar[dr]&&{\begin{smallmatrix}3\\1\;7\\5\;\;\;\;\\6\;\;\;\;\end{smallmatrix}}\ar[ur]\ar[dr]&&{\begin{smallmatrix}3\\1\\5\end{smallmatrix}}\ar[ur]\ar[dr]&&{\begin{smallmatrix}8\\3\\1\end{smallmatrix}}\ar[ur]\ar[dr]&&{\begin{smallmatrix}4\end{smallmatrix}}\\&{\begin{smallmatrix}P_{10}[1]\end{smallmatrix}}\ar[ur]\ar[dr]&&{\begin{smallmatrix}10\\2\;11\\3\;\;\;\;\;\\7\;\;\;\;\;\end{smallmatrix}}\ar[ur]\ar[dr]&&{\begin{smallmatrix}10\\2\\3\end{smallmatrix}}\ar[ur]\ar[dr]&&{\begin{smallmatrix}8\end{smallmatrix}}\ar[ur]\ar[dr]&&{\begin{smallmatrix}7\end{smallmatrix}}\ar[ur]\ar[dr]&&{\begin{smallmatrix}3\\1\\5\\6\end{smallmatrix}}\ar[ur]\ar[dr]&&{\begin{smallmatrix}8\\3\\1\\5\end{smallmatrix}}\ar[ur]\ar[dr]&&{\begin{smallmatrix}P_5[1]\end{smallmatrix}}\ar[ur]\ar[dr]\\{\begin{smallmatrix}P_{11}[1]\end{smallmatrix}}\ar[ur]\ar[dr]&&{\begin{smallmatrix}11\end{smallmatrix}}\ar[ur]\ar[dr]&&{\begin{smallmatrix}10\\2\\3\\7\end{smallmatrix}}\ar[ur]\ar[dr]&&{\begin{smallmatrix}P_{7}[1]\end{smallmatrix}}\ar[ur]\ar[dr]&&{\begin{smallmatrix}7\\8\end{smallmatrix}}\ar[ur]\ar[dr]&&{\begin{smallmatrix}P_{8}[1]\end{smallmatrix}}\ar[ur]\ar[dr]&&{\begin{smallmatrix}8\\3\\1\\5\\6\end{smallmatrix}}\ar[ur]\ar[dr]&&{\begin{smallmatrix}P_{6}[1]\end{smallmatrix}}\ar[ur]\ar[dr]&&{\begin{smallmatrix}6\end{smallmatrix}}\\&{\begin{smallmatrix}P_{x_{2}}[1]\end{smallmatrix}}\ar[ur]&&{\begin{smallmatrix}P_{x_{1}}[1]\end{smallmatrix}}\ar[ur]&&{\begin{smallmatrix}P_{x_{14}}[1]\end{smallmatrix}}\ar[ur]&&{\begin{smallmatrix}P_{x_{13}}[1]\end{smallmatrix}}\ar[ur]&&{\begin{smallmatrix}P_{x_{12}}[1]\end{smallmatrix}}\ar[ur]&&{\begin{smallmatrix}P_{x_{11}}[1]\end{smallmatrix}}\ar[ur]&&{\begin{smallmatrix}P_{x_{10}}[1]\end{smallmatrix}}\ar[ur]&&{\begin{smallmatrix}P_{x_{9}}[1]\end{smallmatrix}}\ar[ur]}$\end{minipage}}}%
  \end{picture}%
\endgroup%
}}

 \caption{AR quiver of the category $\mathcal C_f$ arising from $Q$}
   \label{fig:ARquiver}
\end{figure}
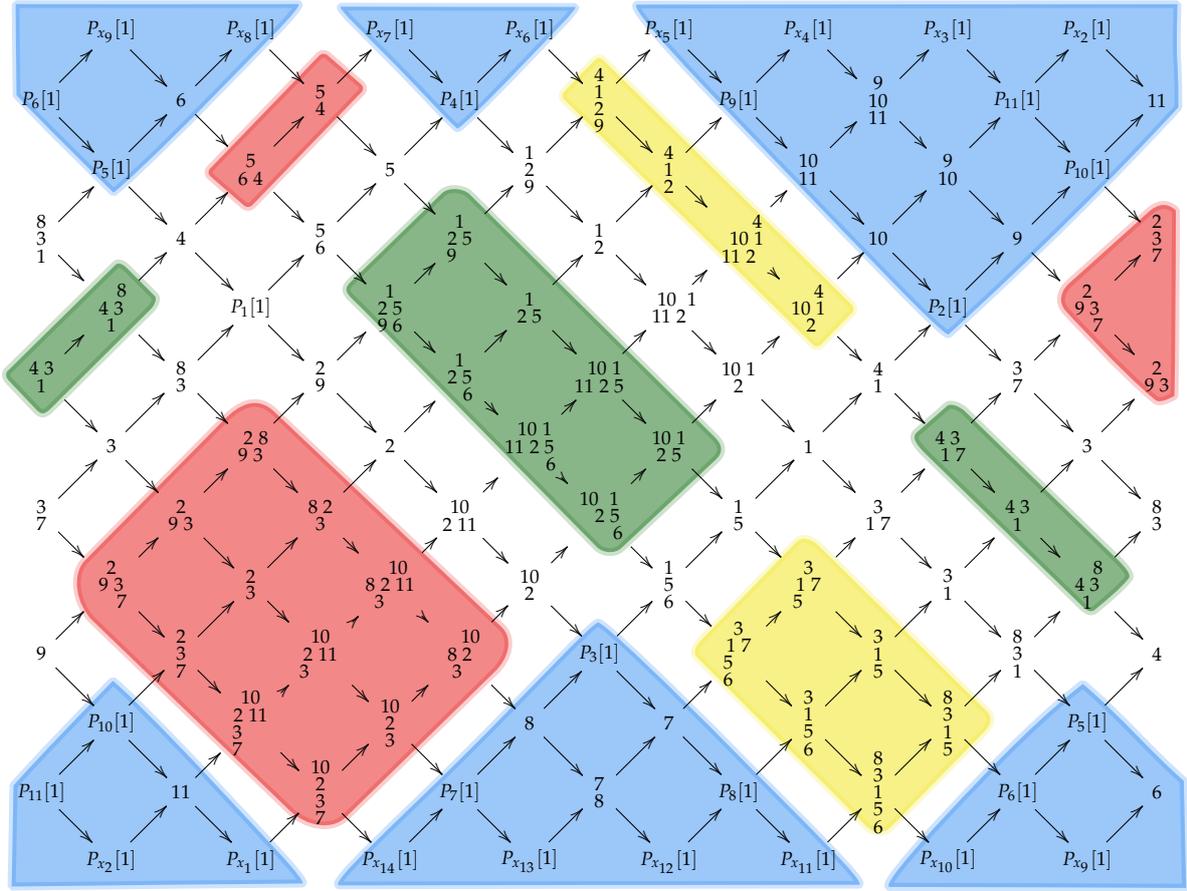

\end{ex}

%%%%%%%%%%%%

\section{Mutating friezes}\label{sec:mutation}

Assume now that our cluster tilting object $T$ in $\CC$ is of the form $T = \bigoplus_{i=1}^n T_i$, where the $T_i$ are mutually non-isomorphic indecomposable objects. Mutating $T$ at $T_i$ for some $1 \leq i \leq n$ yields a new cluster tilting object $T' = T/T_i \oplus T'_i$, to which we can associate a new 
frieze $F(T')$. 
In terms of the frieze, we can think of this mutation as a mutation at an entry of value $1$, 
namely the one sitting in the position of the indecomposable object $T_i[1]$.

We describe how, using graphic calculus, we can obtain each entry of the frieze $F(T')$ independently and directly from the frieze $F(T)$, thus effectively introducing the concept of mutations of friezes at entries of value $1$ that do not lie in one of the two constant rows of $1$s at bounding the frieze pattern. 

We are able to give an explicit formula 
of how each entry in the frieze $F(T)$ changes under mutation at the entry corresponding to $T_i$, 
see Theorem~\ref{thm01} below. 
We observe that each frieze can be divided into four separate regions, relative to the entry of value $1$ at which we want to mutate. Each of these regions gets affected differently by mutation. The formula of the 
theorem relies solely on the shape of the frieze and the entry at which we mutate. 
It determines how each entry of the frieze individually changes under mutation.

In Section~\ref{sec:mutating} we will describe the four separate regions in 
and introduce the necessary notation before stating the theorem. 

%%%%%%%%%%%%
\subsection{Frieze category}% and regions} 

We extend $\ind\mathcal C$ by adding an indecomposable for each boundary segment of 
the polygon and denote the resulting category 
by $\mathcal{C}_f$. Then $\mathcal C_f$ is the Frobenius category of maximal CM-modules 
categorifying 
the cluster algebra structure of the coordinate 
ring of the (affine cone of the) Grassmannian Gr(2,n) as studied in \cite{DeLuo} 
and for general Grassmannians in \cite{JKS,BKM}. 
The stable category of $\mathcal{C}_f$ is equivalent to $\mathcal{C}$. 
We then extend the definition of $\rho_T$ to $\mathcal{C}_f$ by setting 
\[
\rho_T(M)=1 \quad\mbox{if $M$ corresponds to a boundary segment.}
\]
This agrees with the extension of the cluster character to Frobenius category 
given by Fu and Keller, cf. Theorem~\cite[Theorem 3.3]{FuKeller}.

%%%%%%
\section{Mutating friezes}\label{sec:mutating}
%%%%%%

The goal of this section is to describe the effect of the flip of a diagonal or equivalently the mutation at 
an indecomposable projective on the associated frieze. We give a formula for computing the effect 
of the mutation 
using the specialised Caldero Chapoton map. 
Let $\mathcal{T}$ be a triangulation of a polygon with associated quiver $Q$ 
(see Section~\ref{sec:regions}). 
The quiver $Q$ looks as in Figure~\ref{fig:quiver-divided}, where the subquivers $Q_b$, $Q_c$, 
$Q_d$, $Q_e$ may 
be empty. Let 
$T=\oplus_{x\in T} P_x$ and  
$B=\End_{\mathcal C} T$ be the associated cluster-tilted algebra. 

\begin{figure}
 \includegraphics[width=45mm]{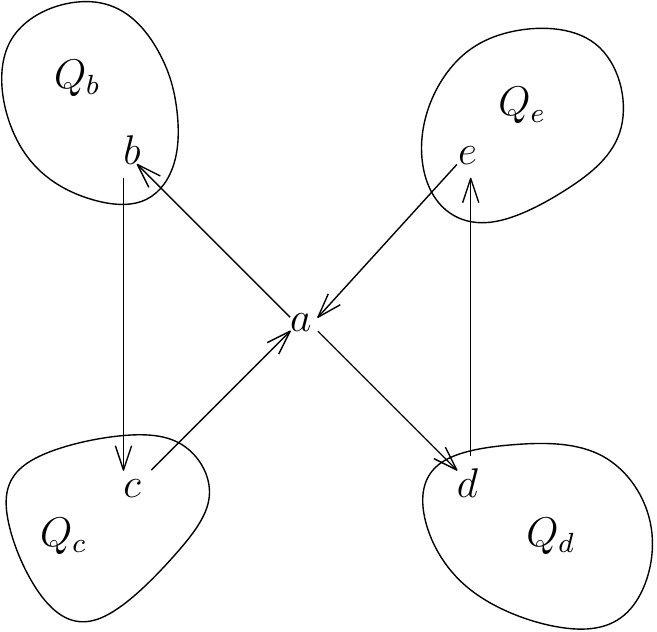}
 \caption{Quiver after flipping diagonal $a$}
 \label{fig:flipped-quiver}
\end{figure}

Take $a\in \mathcal{T}$ and let $\mathcal{T}'=\mu_a(\mathcal{T})$ be the triangulation obtained from flipping $a$, 
with quiver $Q'=\mu_a(Q)$ (Figure~\ref{fig:flipped-quiver}). 
%\[
%\]

Let $B'$ be the algebra obtained through this, it is the cluster-tilted algebra for 
$T'=\oplus_{x\in \mathcal{T}'}P_x$. If $M$ is an indecomposable $B$-module, we write $M'$ for 
$\mu_{a} (M)$ in the sense of \cite{DWZ1}. If $M$ is an indecomposable $B$-module, the {\em entry of $M$} in the frieze $F({T})$ is the 
entry at the position of $M$ in the frieze.

\begin{defi}
Let $\mathcal{T}$ be a triangulation of a polygon, $a\in \mathcal{T}$ and $M$ an indecomposable 
object of $\mathcal{C}_f$. Then we define the {\em frieze difference (w.r.t. mutation at $a$)}  
$\delta_a: \ind \mathcal{C}_f \to \mathbb{Z}$ by 
 \[
 \delta_a(M)=\rho_{\mathcal{T}}(M)-\rho_{\mathcal{T}'}(M') \in \mathbb{Z} 
 \]
\end{defi}

In Section~\ref{ssec:mutating-regions} we first describe the effect mutation has on 
the regions in the frieze. This gives us all the necessary tools to compute the frieze difference 
$\delta_a$ (Section~\ref{ssec:computing-delta}). 

%%%%%%
\subsection{Mutation of regions}\label{ssec:mutating-regions}
%%%%%%
$\ $

Here we describe how mutation affects the regions (Section~\ref{ssec:regions}) of the frieze $F({T})$. 
Let $\mathcal{T},a,B$ and $\mathcal{T}',B'$ be as above. When mutating at $a$, the change in support of 
the indecomposable modules can be described explicitly in terms of the local quiver around 
$a$. This is what we will do here. We first describe the regions in the AR quiver of $\mathcal{C}_f$ 
for $B'$.

If $x$ is a diagonal or a boundary segment, we write 
\[
\mathcal X'=\{M \in \mbox{ind}\,B'\mid \Hom(P_x,M)\ne 0\}
\]
for the indecomposable modules supported on $x$.

After mutating $a$, the regions in the AR quiver are still determined by the projective indecomposables corresponding 
to the framing diagonals (or edges) $b,c,d,e$. The relative positions of $a,b,c,d$ and $e$ have changed, however it follows from \cite{DWZ1} that 
except for vertex $a$ the support of an indecomposable module at all other vertices remains the same.  
Therefore, the regions 
are now described as follows: 

$$\mathcal{B}'\cap \mathcal{E}' = \{ M \in \text{ind}\,B' \mid M \text{ is supported on } e\to a\to b \}$$

$$\mathcal{C}'\cap \mathcal{D}' = \{ M \in \text{ind}\,B' \mid M \text{ is supported on } c\to a \to d \}$$

$$\mathcal{B}'\cap \mathcal{C}' = \{ M \in \text{ind}\,B' \mid M \text{ is supported on } b\to   c \}$$

$$\mathcal{D}'\cap \mathcal{E}'= \{ M \in \text{ind}\,B' \mid M \text{ is supported on } d\to   e \}$$

$$\mathcal{B}'\cap \mathcal{D}' = \{ M \in \text{ind}\,B' \mid M \text{ is supported on } b\leftarrow  a  \to d \}$$

$$\mathcal{C}'\cap \mathcal{E}' = \{ M \in \text{ind}\,B' \mid M \text{ is supported on } c \to a  \leftarrow e \}$$

Under the mutation at $a$, if a module $M$ lies on one of the rays $\mathfrak{b}_a$, $\mathfrak{d}_a$ $\mathfrak{c}^a$ and $\mathfrak{e}^a$ 
then $M'$ is obtained from $M$ by removing support at vertex $a$.  The modules lying on the remaining four rays gain support at vertex $a$ after the mutation.

%%%%%%
\subsection{Mutation of frieze}\label{ssec:computing-delta}
%%%%%%
$\ $

We next present the main result of this section, the effect of flip on the generalized 
Caldero Chapoton map, i.e. the description of the frieze difference 
$\delta_a$. 
We begin by introducing the necessary notation. 

Depending on the position of 
an indecomposable object $M$ we define several projection maps sending 
$M$ to objects on the eight rays from Section~\ref{ssec:diagonal-rays}.

Let $M\in\text{ind}\,B$, and let $\mathfrak{i}$ be one of the sectional paths defined in 
section~\ref{ssec:diagonal-rays}.  Suppose $M \not\in\mathfrak{i}$, then we denote by 
$M_i$ a module on $\mathfrak{i}$ if there exists a sectional path $M_i \to \dots \to M$ or 
$M\to \dots \to M_i$ in $\mathcal{C}_f$, otherwise we let $M_i=0$.  If $M \in\mathfrak{i}$ 
then we let $M_i=M$.  
In the case when it is well-defined, we call $M_i$ the \emph{projection} of $M$ onto 
the path $\mathfrak{i}$.  

It will be convenient to write these projections in a uniform way.

\begin{figure}
\includegraphics[height=5.5cm]{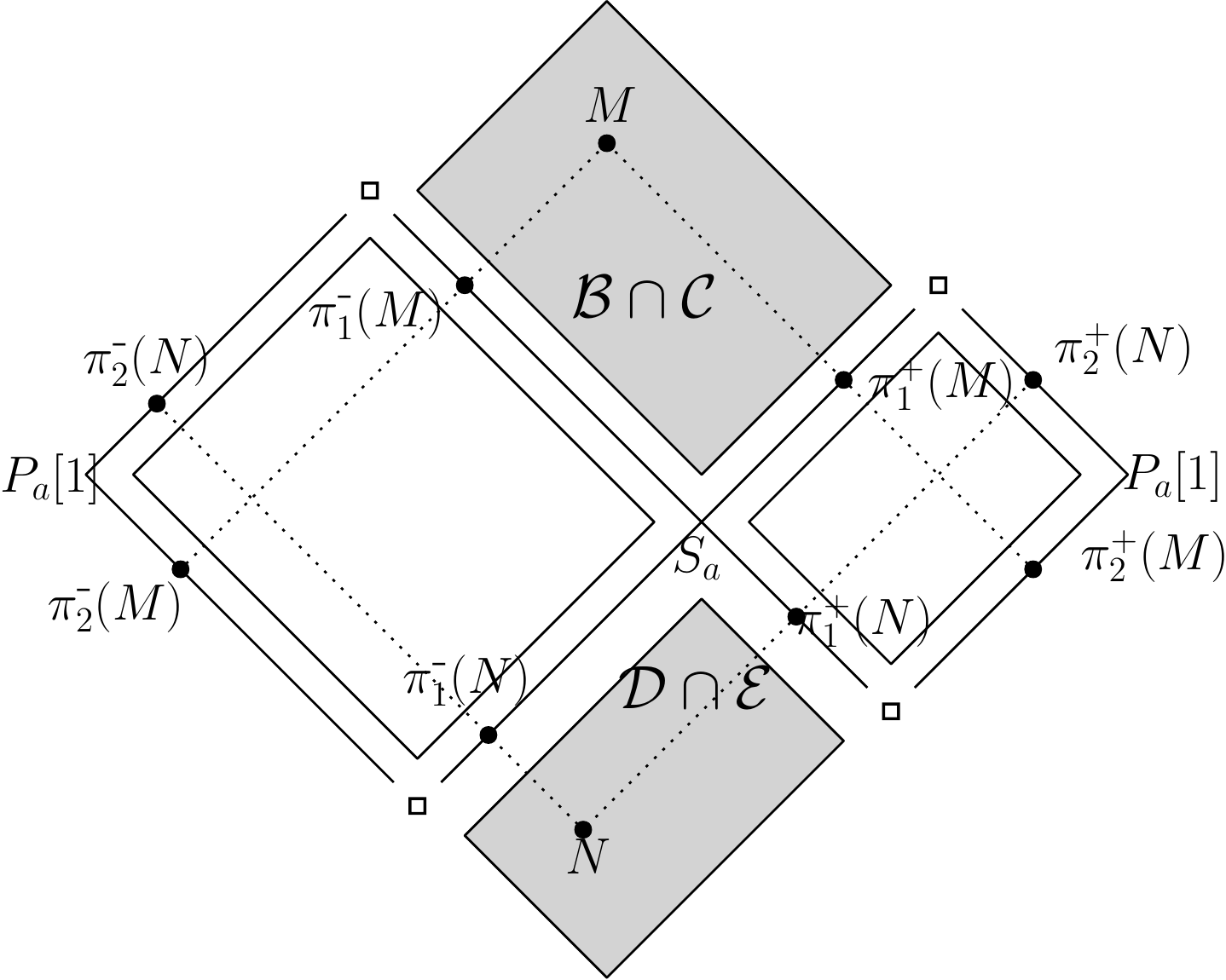}
\caption{Projections for $\mathcal B\cap \mathcal C$, $\mathcal D\cap \mathcal E$}\label{fig:projections-case-a}
\end{figure}

\begin{defi}[Projections] 
If $(x,y)$ is one of the pairs 
$\{(b,c),(d,e),(b,e),(c,d)\}$, the region $\mathcal X\cap \mathcal Y$ 
has two paths along its boundary and 
two paths further backwards or forwards met along the two sectional paths through any vertex $M$ of 
$\mathcal X\cap \mathcal Y$.
We call the backwards projection onto the first path 
$\pi_1^-(M)$ and the projection onto the second path 
$\pi_2^-(M)$. The forwards 
projection onto the first path is denoted by $\pi_1^+(M)$ and the one onto the second path $\pi_2^+(M)$. 
\end{defi}
Figure~\ref{fig:projections-case-a} illustrates these projections in the case 
$(x,y)\in \{(b,c),(d,e)\}$.

\smallskip

The remaining two regions will be treated together with the surrounding paths. 

\begin{defi}
The {\em closure of $\mathcal C\cap \mathcal E$} is the 
Hom-hammock 
$$
\overline{\mathcal C\cap \mathcal E}=\ind(\Hom_{\mathcal C_f}(P_a[1],-)\cap \Hom_{\mathcal C_f}(-,S_a))
$$ 
in $\mathcal C_f$ starting at $P_a[1]$ and ending at $S_a$. 
Similarly, 
the {\em closure of $\mathcal B\cap \mathcal D$} is the Hom-hammock 
$$
\overline{\mathcal B\cap \mathcal D}=\ind(\Hom_{\mathcal C_f}(S_a,-)\cap \Hom_{\mathcal C_f}(-,P_a[1]))
$$
in $\mathcal C_f$ starting at $S_a$ and ending at $P_a[1]$. 
For $(x,y)\in \{(c,e),(b,d)\}$, the {\em boundary of $\overline{\mathcal X\cap \mathcal Y}$} 
(or of 
$\mathcal X\cap \mathcal Y$) is 
$\overline{\mathcal X\cap \mathcal Y}\setminus(\mathcal X\cap \mathcal Y)$.  
\end{defi}

Note that $\overline{\mathcal C\cap \mathcal E}$ is the union of 
$\mathcal C\cap \mathcal E$ with the surrounding rays and the shifted projectives 
$\{P_b[1],P_d[1]\}$. Analogously, 
$\overline{\mathcal B\cap \mathcal D}$ contains $\{P_c[1],P_e[1]\}$. 

\begin{figure}
\includegraphics[height=5cm]{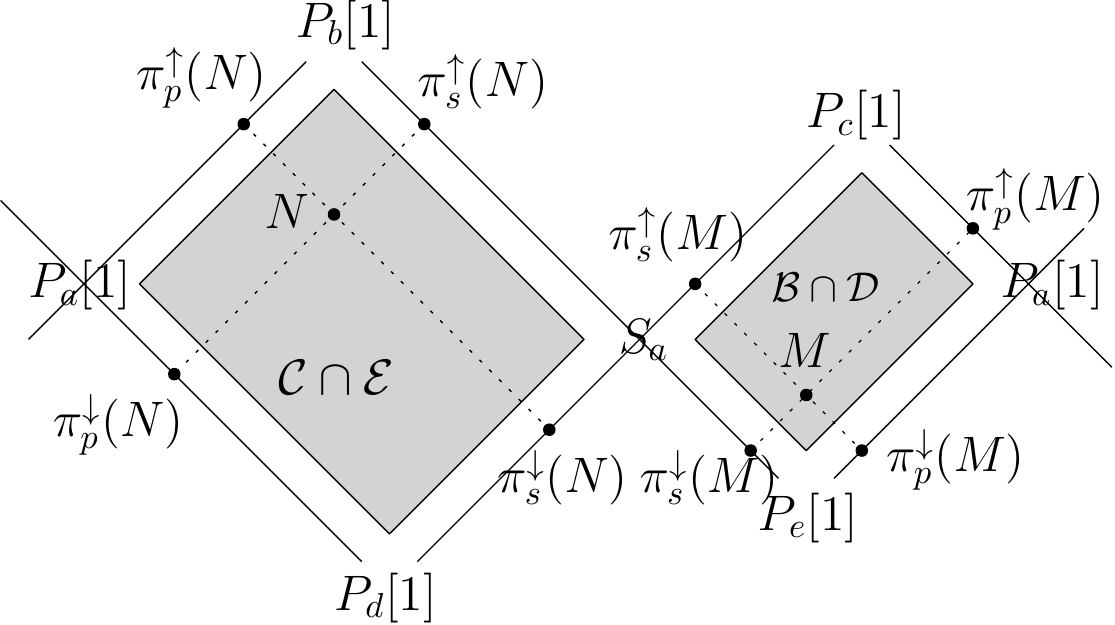}
\caption{}\label{fig:projections-bar}
\end{figure}
\begin{defi}[Projections, continued] \label{def:projections-bar}
If $M$ is a vertex of one of the two closures, 
we define 
four projections for $M$ onto the four different ``edges'' of the boundary of 
its region: 
We denote the 
projections onto the paths starting or ending next to $P_a[1]$ by 
$\pi_{p}^{\uparrow}$, $\pi_{p}^{\downarrow}$ 
and the projections onto the paths starting or ending next to $S_a$ 
by $\pi_{s}^{\uparrow}$ and $\pi_{s}^{\downarrow}$ respectively. 
We choose the upwards arrow to refer to the paths ending/starting near $P_b[1]$ or 
$P_c[1]$ and the downwards arrow to refer to paths ending/starting near 
$P_d[1]$ or $P_e[1]$.  
See Figure~\ref{fig:projections-bar}. 
\end{defi}

\begin{rem}
The statement of Theorem~\ref{thm01} is independent of the 
choice of  $\uparrow$ (paths near $P_b[1]$ or $P_c[1]$) and 
$\downarrow$ in Definition~\ref{def:projections-bar} 
as the formula is symmetric in these expressions. 
\end{rem}

\begin{ex}\label{ex-projections}
If $M\in \mathfrak e$, we have 
$\pi_p^{\uparrow}(M)=M$, $\pi_s^{\uparrow}(M)=P_b[1]$, $\pi_p^{\downarrow}(M)=P_a[1]$ 
and $\pi_s^{\downarrow}(M)=M_{\mathfrak e^a}$. 

For $S_a$  we have 
$\pi_s^{\uparrow}(S_a)=\pi_s^{\downarrow}(S_a)=S_a$ whereas the two modules 
$\pi_p^{\uparrow}(S_a)$ and $\pi_p^{\downarrow}(S_a)$ are  $\{P_b[1],P_d[1]\}$ or 
$\{P_c[1],P_e[1]\}$ depending on whether $S_a$ is viewed as an element of 
$\overline{\mathcal C\cap \mathcal E}$ or of $\overline{\mathcal B\cap \mathcal D}$. 

For $P_a[1]$, we have 
$\pi_p^{\uparrow}(P_a[1])=\pi_p^{\downarrow}(P_a[1])=P_a[1]$ whereas the two modules 
$\pi_s^{\uparrow}(P_a[1])$ and $\pi_s^{\downarrow}(P_a[1])$ are $\{P_b[1],P_d[1]\}$ or 
$\{P_c[1],P_e[1]\}$
These four shifted projectives evaluate to $1$ under $s$, and so in Theorem~\ref{thm01}, 
this ambiguity does not play a role. 
\end{ex}

With this notation we are ready to state the theorem, proved in~\cite[Section 6]{BFGST17}. 
 
\begin{theorem}\label{thm01}
Consider a frieze associated to a triangulation of a polygon. Let $a$ be a diagonal in the triangulation.
Consider the mutation of the frieze at $a$. Then the frieze difference $\delta_a(M)$ at the point corresponding to the indecomposable object $M$ in the associated Frobenius category $\mathcal{C}_f$ is given by:\\\\
%Let $M$ be an indecomposable object of $\mathcal{C}_f$. Then $\delta_a(M)$ is given by: \\
If $M\in (\mathcal B\cap \mathcal C) \cup (\mathcal D\cap \mathcal E)$ then
\[
 \delta_a(M) = (s(\pi_1^+(M))-s(\pi_2^+(M)))\,(s(\pi_{1}^-(M))-s(\pi_{2}^-(M));
\]
If $M\in (\mathcal B\cap \mathcal E) \cup (\mathcal C\cap \mathcal D)$ then
\[
 \delta_a(M) = -(s(\pi_2^+(M))-2s(\pi_1^+(M)))\,(s(\pi_{2}^-(M))-2s(\pi_{1}^-(M)) ;
\]
If $M\in \overline{\mathcal C\cap \mathcal E}\cup\overline{\mathcal B\cap \mathcal D}$ then
\[
 \delta_a(M) = s(\pi_{s}^{\downarrow}(M))s(\pi_{p}^{\downarrow}(M)) + s(\pi_{s}^{\uparrow}(M))s(\pi_{p}^{\uparrow}(M))  - 3\,s(\pi_p^{\downarrow}(M))s(\pi_p^{\uparrow}(M)); 
\]
If $M\in \mathcal F$ then
\[
 \delta_a(M) =0.
\]

\end{theorem}

Note, that given a frieze 
and an indecomposable $M$ in one of the six regions $\mathcal X\cap \mathcal Y$, 
it is easy to locate the entries required to compute the frieze difference $\delta_a(M)$. 
We simply need to find projections onto the appropriate rays in the frieze. 
In this way, we do not need to know the precise shape of the modules appearing in the formulas of Theorem~\ref{thm01}.

\begin{ex}
Let $\mathcal{C}_f$ be the category given in Example~\ref{ex:quiver-triangul-flip}.  We consider three possibilities for $M$ below.   

If $M = \begin{smallmatrix} \;\;\;\;\;\;4\\\;\;10\;1\\11\;2\end{smallmatrix}$ then we know by Figure~\ref{fig:frieze-entries} that $s(M)= 11$ and $s(M')= 8$.   On the other hand, we see from Figure~\ref{fig:ARquiver} that $M \in \mathcal{B}\cap\mathcal{C}$.  Theorem~\ref{thm01} implies that 

\begin{align*}
\delta_a(M)=s(M)-s(M')&=(s(M_{b_a})-s(M_b))(s(M_{c^a})-s(M_c))\\
& = (s(\begin{smallmatrix} 4\\1\end{smallmatrix})-s(\begin{smallmatrix} 4\end{smallmatrix}))(s(\begin{smallmatrix}\;\;10\;1\\11\;2\end{smallmatrix})-s(\begin{smallmatrix} 10\\11\;2 \end{smallmatrix}))\\
&=(3-2)(8-5)=3.
\end{align*}

Similarly, if $M = \begin{smallmatrix}8\;2\\3\end{smallmatrix}$, then $M \in \mathcal{C}\cap\mathcal{D}$ 
with $s(M)=5$ and $s(M')=7$.  The same theorem implies that

\begin{align*}
\delta_a(M)=s(M)-s(M') &= -(s(M_{c^a})-2s(M_{c}))\,(s(M_{d^a})-2s(M_d))\\
&=-(s(\begin{smallmatrix}1\\2\end{smallmatrix})- 2s(\begin{smallmatrix}2\end{smallmatrix})) 
(s(\begin{smallmatrix}8\\3\\1\end{smallmatrix})-2s(\begin{smallmatrix}8\\3\end{smallmatrix}))\\
&=-(3-4)(4-6)=-2.
\end{align*}

Finally, if $M=\begin{smallmatrix}10\;1\;\\\;\;2\;5\\\;\;\;\;\;\;6\end{smallmatrix}$, then $M \in \mathcal{C}\cap\mathcal{E}$.  We also know that $s(M)=s(M')=11$.  By the third formula in Theorem~\ref{thm01}, we have 

\begin{align*}
\delta_a(M)=s(M)-s(M')&=s(M_{e^a})s(M_c) + s(M_{c^a})s(M_e) - 3 s(M_e) s(M_c)\\
&= s(\begin{smallmatrix}1\\5\\6\end{smallmatrix})s(\begin{smallmatrix}10\\2\end{smallmatrix}) + s(\begin{smallmatrix}10\;1\\2\end{smallmatrix}) s(\begin{smallmatrix}5\\6\end{smallmatrix}) -3 s(\begin{smallmatrix}5\\6\end{smallmatrix}) s(\begin{smallmatrix}10\\2\end{smallmatrix}) \\
&= 4\cdot 3 + 5\cdot 3 - 3 \cdot 3 \cdot 3 = 0.
\end{align*}
\end{ex}

%%%%%%%%%%%%

\section*{Acknowledgements}

We thank 
AWM for encouraging us to write this summary and giving us opportunity to continue 
this work. \\
EF, KS and GT received support from the AWM Advance grant to attend the symposium.\\
KB was supported by the FWF grant W1230. \\
KS was supported by NSF Postdoctoral Fellowship MSPRF - 1502881.

%%%%%%%%bibtex 
%\bibliographystyle{alpha}
%\bibliography{biblio}

\newcommand{\etalchar}[1]{$^{#1}$}
\def\cprime{$'$}

\end{document}